\documentclass{gtart_h}

\def\ifplaintex{\expandafter\ifx\csname documentclass\endcsname\relax}

\def\gtp{{\mathsurround=0pt\it $\cal G\mskip-2mu$eometry \&\ 
$\cal T\!\!$opology $\cal P\!$ublications}}  

\def\recd{{\small Received:\qua\receiveddate\ifx\reviseddate\relax
\else\qquad Revised:\qua\reviseddate\fi\par}} 

\def\lognumber#1{\def\thelognumber{#1}}
\def\volumenumber#1{\def\thevolumenumber{#1}}
\def\volumeyear#1{\def\thevolumeyear{#1}}
\def\papernumber#1{\def\thepapernumber{#1}}
\def\pagenumbers#1#2{\def\startpage{#1}\def\finishpage{#2}}
\def\published#1{\def\publishdate{#1}}

\def\received#1{\def\receiveddate{#1}}

\def\accepted#1{\def\accepteddate{#1}}

\def\asciiurl#1{\def\theasciiurl{#1}}

\long\def\asciiabstract#1{\long\def\theasciiabstract{#1}}

\let\\\par\let\thelognumber\relax\let\thevolumenumber\relax
\let\thepapernumber\relax\let\thevolumeyear\relax\let\startpage\relax
\let\finishpage\relax\let\publishdate\relax\let\receiveddate\relax
\let\reviseddate\relax\let\accepteddate\relax\let\theasciititle\relax
\let\theasciiauthors\relax
\let\theasciiabstract\relax

\let\theasciiemail\relax
\let\theasciiurl\relax

\ifplaintex
\font\logobig=cmssbx10 scaled 3836
\font\logomed=cmssbx10 scaled 2557
\else
\font\logobig=cmssbx10 scaled 4200
\font\logomed=cmssbx10 scaled 2800
\fi

\long\def\makeagttitle{   
\count0=\startpage
\agt\hfill      
\hbox to 45truept{\vbox to 0pt{\vglue -13truept{\logomed A\kern -.37em{\logobig 
T}\kern -.38em G}\vss}\hss}
\break
{\small Volume \thevolumenumber\ (\thevolumeyear)
\startpage--\finishpage\nl
Published: \publishdate}

\vglue .25truein

{\parskip=0pt\leftskip 0pt plus
1fil\def\\{\par\smallskip}{\Large\bf\thetitle}\par\medskip} \vglue
0.05truein

{\parskip=0pt\leftskip 0pt plus 1fil\def\\{\par}{\sc\theauthors}
\par\medskip}%
 
\vglue 0.03truein 

{\small\leftskip 25truept\rightskip 25truept{\bf Abstract}\stdspace\theabstract

{\bf AMS Classification}\stdspace\theprimaryclass
\ifx\thesecondaryclass\relax\else; \thesecondaryclass\fi\par
{\bf Keywords}\stdspace \thekeywords\par}\vglue 7truept

}   

\ifplaintex
\font\phead=cmsl9 scaled 950
\font\pnum=cmbx10 scaled 913
\font\pfoot=cmsl9 scaled 950
\headline{\vbox to 0pt{\vskip -4.5mm\line{\small\phead\ifnum
\count0=\startpage ISSN 1472-2739 (on-line) 1472-2747 (printed)
\hfill {\pnum\folio}\else\ifodd\count0\def\\{ }%
\ifx\theshorttitle\relax\thetitle\else\theshorttitle\fi\hfill{\pnum\folio}
\else\def\\{ and }{\pnum\folio}\hfill\ifx\theshortauthors\relax\theauthors
\else\theshortauthors\fi\fi\fi}\vss}}
\footline{\vbox to 0pt{\vglue 0mm\line{\small\pfoot\ifnum\count0=\startpage
\copyright\ \gtp\hfill\else
\agt, Volume \thevolumenumber\ (\thevolumeyear)\hfill\fi}\vss}}
\else
\font=cmsl9 scaled 1050
\font\lnum=cmbx10 
\font=cmsl9 scaled 1050
\makeatletter
\def\@oddhead{{\small\lhead\ifnum\count0=\startpage ISSN 1472-2739 
(on-line) 1472-2747 (printed)\hfill {\lnum\number\count0}\else\ifodd\count0
\def\\{ }\ifx\theshorttitle\relax \thetitle \else\theshorttitle\fi\hfill
{\lnum\number\count0}\else\def\\{ and }{\lnum\number\count0}
\hfill\ifx\theshortauthors\relax 
\theauthors\else\theshortauthors\fi\fi\fi}}\def\@evenhead{\@oddhead}
\def\@oddfoot{\small\lfoot\ifnum\count0=\startpage\copyright\ \gtp\hfill\else
\agt, Volume \thevolumenumber\ (\thevolumeyear)\hfill\fi}
\def\@evenfoot{\@oddfoot}
\makeatother
\fi
\let\maketitlepage\makeagttitle
\let\makeshorttitle\maketitlepage
\let\maketitle\maketitlepage

\newwrite\gtoutfile
\long\gdef\makeheadfile{  
{\def\\{, }\def\s{ }
\immediate\openout\gtoutfile head.xxx
\immediate\write\gtoutfile{Proxy-for: \ifx\theasciiauthors\relax
\theauthors\else\theasciiauthors\fi\s<\ifx\theasciiemail\relax\theemail\else\theasciiemail\fi>}
\immediate\write\gtoutfile{\noexpand\\}
\immediate\write\gtoutfile{Authors: \ifx\theasciiauthors\relax
\theauthors\else\theasciiauthors\fi}
{\def\\{ }\immediate\write\gtoutfile{Title: \ifx\theasciititle\relax
\thetitle\else\theasciititle\fi}}
\immediate\write\gtoutfile{Subj-class: GT or SG, GR etc}
\immediate\write\gtoutfile{MSC-class: \theprimaryclass\ifx\thesecondaryclass\relax\else, \thesecondaryclass\fi}
\immediate\write\gtoutfile{Journal-ref: Algebr. Geom. Topol. \thevolumenumber\s
(\thevolumeyear) \startpage-\finishpage}
\immediate\write\gtoutfile{Comments: Published by Algebraic and
Geometric Topology at}
\immediate\write\gtoutfile{\s\s\s  http://www.maths.warwick.ac.uk/agt/AGTVol\thevolumenumber/agt-\thevolumenumber-\thepapernumber.abs.html}
\immediate\write\gtoutfile{\noexpand\\}
\immediate\write\gtoutfile{}
\ifx\theasciiabstract\relax
\immediate\write\gtoutfile{\theabstract}\else
\immediate\write\gtoutfile{\theasciiabstract}\fi
\immediate\write\gtoutfile{}
\immediate\write\gtoutfile{\noexpand\\}
\immediate\write\gtoutfile{}
\immediate\closeout\gtoutfile}}  

\def\maketitlepage{\makeagttitle\makeheadfile}
\let\makeshorttitle\maketitlepage
\let\maketitle\maketitlepage

\lognumber{45}
\volumenumber{5}
\volumeyear{2005}
\papernumber{45}
\pagenumbers{1111}{1139}
\received{6 May 2005} 
\accepted{15 August 2005}
\published{11 September 2005}

\usepackage{amsmath,amssymb,graphicx,epsfig,verbatim,color}
\usepackage[all]{xy}

\newcommand{\Z}{\ensuremath{\mathbb{Z}}}

\newcommand{\Uh}{\ensuremath{U_{h}(\mathfrak{g}) } }

\newcommand{\s}{\ensuremath{\operatorname{s}}}
\newcommand{\g}{\ensuremath{\mathfrak{g}}}

\newcommand{\gl}{\ensuremath{\mathfrak{gl}(m|n)}}

\newcommand{\gto}{\ensuremath{\mathfrak{gl}(2|1)}}
\newcommand{\C}{\ensuremath{\mathbb{C}} }
\newcommand{\Q}{\ensuremath{\mathbb{Q}} }
\newcommand{\N}{\ensuremath{\mathbb{N}} }
\newcommand{\p}[1]{\ensuremath{\bar {#1}}}

\newcommand{\A}{\ensuremath{\mathcal{A}}}

\newcommand{\unit}{\ensuremath{\mathbb{I}}}
\newcommand{\W}{\ensuremath{\mathcal{W}}}
\newcommand{\V}{\ensuremath{\mathcal{V}}}
\newcommand{\chord}{\ensuremath{\mathcal{D}^{c}}}
\newcommand{\chordm}{\ensuremath{\mathcal{D}_{m}^{c}}}
\newcommand{\R}{\ensuremath{\mathbb{R}} }
\newcommand{\cat}{\ensuremath{\mathcal{C}}}
\newcommand{\rib}{\ensuremath{\mathcal{R}}}
\newcommand{\sym}{\ensuremath{\mathcal{S}}}

\newcommand{\Hom}{\operatorname{Hom}}
\newcommand{\End}{\operatorname{End}}
\newcommand{\Atangle}{\ensuremath{A}}

\newcommand{\kzt}{A_{\g,t}}

\newcommand{\lgws}{\ensuremath{W_{LG}}}
\newcommand{\twist}{\theta}
\newcommand{\Valpha}{V_{\alpha}}
\newcommand{\Valphahat}{\widetilde{V}_{\alpha}}
\newcommand{\cT}{(T^{\otimes 2})'}
\newcommand{\KontFunc}{\zeta}

\newcommand{\stan}{M}

\newcommand{\T}[1]{\widehat{#1}}
\newcommand{\Vhat}{\ensuremath{\T{V}}}

\newcommand{\flip}{\operatorname{\tau}}
\newcommand{\slmn}{\ensuremath{\mathfrak{sl}(m|n)}}
\newcommand{\Vws}{\ensuremath{\Phi}}
\newcommand{\dtoa}{\ensuremath{\mathcal{D}(2,1;\alpha)}}
\newcommand{\Qbar}{\overline{Q}}
\newcommand{\OQ}{\widetilde{Q}}
\newcommand{\OZ}{\widetilde{Z}}

\newcommand{\dto}{\ensuremath{\mathcal{D}(2,1)}}
\newcommand{\dm}{\ensuremath{N}}
\newcommand{\adm}{\ensuremath{M}}
\newcommand{\Qnet}{\ensuremath{Q^{3-tangles}_{\dto,\dm}}}
\newcommand{\Gnet}{\ensuremath{\Phi_{\dto}}}
\newcommand{\ZGad}{\ensuremath{Z_{\dto,\dm}}}

\newcommand{\Etangle}{\ensuremath{E}}

\theoremstyle{definition}
\newtheorem{Df}{Definition}
\newtheorem{remark}[Df]{Remark}

\theoremstyle{plain}
\newtheorem{theorem}[Df]{Theorem}

\newtheorem{prop}[Df]{Proposition}
\newtheorem{lemma}[Df]{Lemma}
\newtheorem{corollary}[Df]{Corollary}

\numberwithin{Df}{section}
\numberwithin{equation}{section}

\begin{document}
\title{The Kontsevich integral and\\quantized Lie superalgebras}

\author{Nathan Geer}
\address{School of Mathematics, Georgia Institute of Technology\\ 
Atlanta, GA 30332-0160, USA}
\email{geer@math.gatech.edu}
\urladdr{www.math.gatech.edu/~geer/}
\asciiurl{www.math.gatech.edu/ geer/}

 \begin{abstract}
Given a finite dimensional representation of a semisimple Lie algebra there are two ways of constructing link invariants:
1) quantum group invariants using the R-matrix,
2) the Kontsevich universal link invariant followed by the Lie algebra based weight system.
Le and Murakami showed that these two link invariants are the same.
These constructions can be generalized to some classes of Lie superalgebras.  In this paper we show that constructions 1) and 2) give the same invariants for the Lie superalgebras of type A-G.   We use this result to investigate the Links-Gould invariant.  We also give a positive answer to a conjecture of Patureau-Mirand's concerning  invariants arising from the Lie superalgebra $\dtoa$.
\end{abstract}
\asciiabstract{%
Given a finite dimensional representation of a semisimple Lie algebra
there are two ways of constructing link invariants: 1) quantum group
invariants using the R-matrix, 2) the Kontsevich universal link
invariant followed by the Lie algebra based weight system.  Le and
Murakami showed that these two link invariants are the same.  These
constructions can be generalized to some classes of Lie superalgebras.
In this paper we show that constructions 1) and 2) give the same
invariants for the Lie superalgebras of type A-G.  We use this result
to investigate the Links-Gould invariant.  We also give a positive
answer to a conjecture of Patureau-Mirand's concerning invariants
arising from the Lie superalgebra D(2,1;alpha).}

  
\maketitle

\section{Introduction}

Given a finite dimensional representation $V$ of a finite dimensional semisimple Lie algebra $\g$ one can construct the following two link invariants: 
\begin{enumerate}
	\item the Reshetikhin-Turaev $\C[[h]]$-valued quantum group invariant $Q_{\g,V}$ which arises from $V$ and the Drinfeld-Jimbo quantization associated to $\g$ \cite{RT, Lin},\label{I:PaperQG}
	\item $W_{\g,V} \circ Z$ where $W_{\g,V}$ is a weight system, constructed by Bar-Natan, and where $Z$ is the Kontsevich integral. \label{I:PaperWS}
\end{enumerate}
These constructions are essentially the same in the following sense.   Lin \cite{Lin} showed that the $m$th coefficient of $Q_{\g,V}$ is a Vassiliev of type $m$.  Moreover, there is a weight system corresponding to  $Q_{\g,V}$ which can be shown to be equal to $W_{\g,V}$.  Conversely, Le and Murakami \cite{LM} showed the invariant $W_{\g,V}\circ Z$ is equal (up to a change of variable and normalization) to  $Q_{\g,V}$.   

For certain classes of Lie superalgebras, there are analogous constructions of the link invariants given in (\ref{I:PaperQG}) and (\ref{I:PaperWS}).   However, in general it is not known that they are essentially the same because, the proof of Le and Murakami uses results, due to Drinfeld, whose proofs are based on properties of Lie algebras which fail for Lie superalgebras.  It was known only for the Lie superalgebra $\mathfrak{gl}(1|1)$.  In this paper we will show that the invariants given in (\ref{I:PaperQG}) and (\ref{I:PaperWS}) are essentially the same for Lie superalgebras of type A-G.    In particular, we will show that 
\begin{equation}
\label{E:PaperIisWSgl}
Q_{\g,V}=W_{\g,V}\circ Z,
\end{equation}
where $\g$ is a Lie superalgebra of type A-G and $V$ is a finite dimensional $\g$-module.
Our proof is based on new quantum group theory results.  In \cite{G1} the author extends the Etingof-Kazhdan theory \cite{EK1} of quantization of Lie bialgebras to Lie superbialgebras and shows that Drinfeld's results hold for the Lie superalgebras of type A-G.

Invariants arising from Lie superalgebras have been studied by many mathematicans.  For example, Vogel \cite{Vogel} proved (at the weight system level) that invariants arising from Lie superalgebras are more powerful than invariants arising from Lie algebras.  From the general linear Lie superalgebra $\mathfrak{gl}(1|1)$ one can recover the Alexander-Conway polynomial (see \cite{FKV}).  In \cite{LG} Links and Gould defined a quantum group invariant $LG$ of (1,1)-tangles arising from quantum superalgebra $U_{q}(\gto)$ and a family of four dimensional simple modules $\Valpha$, $(\alpha \in \C)$.   This invariant is a two variables polynomial.  In \cite{DWLK} De Wit, Links and Kauffman gave an effective way of computing $LG$ using a computer.  Using these calculations they showed it is independent from Jones, HOMFLY and Kauffman polynomials (detecting chirality of some links where these invariants fail).  Recent work on the Links-Gould invariant has been done by Ishii and Kanenobu \cite{Ishii,IshiiKan}.   In \cite{Ishii} Ishii proved that after a variable reduction, the Links-Gould invariant is the Alexander-Conway polynomial.   Finally, Patureau-Mirand \cite{Pat-Mir} modifies the normal construction of quantum invariants to construct a non-zero link invariant from the Lie superalgebra $\dtoa$.

To demonstrate the usefulness of equality (\ref{E:PaperIisWSgl}) we will investigate both the Links-Gould invariant and the invariants defined by Patureau-Mirand.  As mentioned above $LG$ is defined as an invariant of $(1,1)$-tangles.  In \S\ref{SS:LGpower}, we will use equality \eqref{E:PaperIisWSgl} to show that $LG$ is an invariants of knots.    Then, in \S\ref{SS:LGandHOMFLY}, we will show that $LG$ is contained in the colored HOMFLY (that is, every two knots distinguishable by $LG$ can be distinguished by the usual HOMFLY applied to some cabling of these knots).  

Let $V$ be any finite dimensional representation of $\gto$.   As another application of equality (\ref{E:PaperIisWSgl}), we will prove by computation in the Grothendieck ring of representations that $Q_{\gto,V}$ is contained in the cablings of $LG$ and $Q_{\gto,U}$, where $U$ is the defining representation of $\gto$ (see \ref{SS:CablingsLG}).  Finally, we use equality (\ref{E:PaperIisWSgl}) to give a positive answer to Conjecture 4.12 of \cite{Pat-Mir}. This implies that the modified quantum invariant defined by Patureau-Mirand is an invariant which is symmetric in three variables.   This symmetry is not evident from the definition of the quantum invariant.

\subsection*{Acknowledgments}   I would like to thank D. Sinha for reading this paper and giving helpful comments.   I am especially grateful to my PhD advisor A. Vaintrob for  many useful conversations and his insight in to the formulation of these questions.  I would also like to thank P. Etingof for pointing out  \cite{EK1} and T. Le for helpful discussions.   


\section{Preliminaries}


Throughout this paper the ground field is always $\C$.  
In this section we review well known background material and notation that will be used in the following sections.

\subsection{Superspaces and Lie superalgebras}\label{SS:SLS}
A \emph{superspace} is a $\Z_{2}$-graded vector space $V=V_{\p 0}\oplus V_{\p 1}$ over $\C$.  We denote the parity of a homogeneous element $x\in V$ by $\p x\in \Z_{2}$.  We say $x$ is even (odd) if $x\in V_{\p 0}$ (resp. $x\in V_{\p 1}$).  Let $V$ and $W$ be superspaces.    Let $\flip_{V,W}:V\otimes W \rightarrow W \otimes V $ be the linear map given by
\begin{equation}
\label{E:Flip}
 \flip_{V,W}(v\otimes w)=(-1)^{\p v \p w}w\otimes v
\end{equation}
for homogeneous $v \in V$ and $w \in W$. 

A \emph{Lie superalgebra} is a superspace $\g=\g_{\p 0} \oplus \g_{\p 1}$ with a superbracket $[\: , ] :\g^{\otimes 2} \rightarrow \g$ that preserves the  $\Z_{2}$-grading, is super-antisymmetric ($[x,y]=-(-1)^{\p x \p y}[y,x]$), and satisfies the super-Jacobi identity (see \cite{K}).  Throughout, all modules will be $\Z_{2}$-graded modules (module structures which preserve the $\Z_{2}$-grading, see \cite{K}). 
 
 For the purpose of this paper when we say a Lie superalgebras of type A-G we will included $\dtoa$.  All of these Lie superalgebras can be given by generators and relations coming from a Cartan matrix (see \cite{G1,Yam94}).  From Proposition 2.5.3 and 2.5.5 of \cite{Kas} there exists a unique (up to constant factor) non-degenerate supersymmetric invariant bilinear form on each Lie superalgebra of type A-G.  

 \subsection{($k$,$l$)-tangles}\label{SS:klTangles}
  A \emph{(k,l)-tangle} $T$ is a oriented tangles in $\R^{2} \times [0,1]$ such that the boundary $\partial T$ of $T$ satisfies the condition 
$$\partial T = T \cap (\R^{2}\times \{0,1\}) = ([k] \times \{0\} \times \{0\}) \cup ([l]\times \{0\} \times \{1\})$$
where $[n]=\{1,2,\ldots,n\}$.  To each ($k$,$l$)-tangle $T$ we assign two sequences $s(T)$ and $b(T)$ of $+$ and $-$ as follows.  If $k=0$ ($l=0$) let $s(T)=\emptyset$ (resp. $b(T)=\emptyset$).  Otherwise, let $s(T)=(\epsilon_{1},\ldots,\epsilon_{k})$ where $\epsilon_{i}=+ $ if the point $(i,0,0)$ is an endpoint or $\epsilon_{i}=-$ if the point $(i,0,0)$ is an origin.  Similarly, let $b(T)=(\epsilon'_{1},\ldots,\epsilon'_{l})$ where $\epsilon'_{i}=+ $ if the point $(i,0,1)$ is an origin or $\epsilon'_{i}=-$ if the point $(i,0,0)$ is an endpoint.   

  Throughout this paper when we say $T$ is a tangle we will mean that $T$ is a ($k$,$l$)-tangle for some $k$ and $l$.

\subsection{Vassiliev invariants and weight systems}\label{SS:VIandWS}
In this subsection we recall the notions of Vassiliev invariants and weight systems \cite{BN, Kas, KRT97}.  To make a consistent theory of Vassiliev invariants of framed links we restrict to framed links with even framings (see \cite{KRT97}).  

Any numerical link invariant $f$ can be inductively extended to an invariant of singular link according to the rule
\begin{equation*} \scalebox{.6}{ \begin{picture}(0,0)%
\includegraphics{ExtDoublePoint.pstex}%
\end{picture}%
\setlength{\unitlength}{3947sp}%
\begingroup\makeatletter\ifx\SetFigFont\undefined%
\gdef\SetFigFont#1#2#3#4#5{%
  \reset@font\fontsize{#1}{#2pt}%
  \fontfamily{#3}\fontseries{#4}\fontshape{#5}%
  \selectfont}%
\fi\endgroup%
\begin{picture}(3675,474)(3151,-2023)
\put(3151,-1861){\makebox(0,0)[lb]{\smash{{\SetFigFont{17}{20.4}{\familydefault}{\mddefault}{\updefault}{\color[rgb]{0,0,0}$f\Big($}%
}}}}
\put(5926,-1861){\makebox(0,0)[lb]{\smash{{\SetFigFont{17}{20.4}{\familydefault}{\mddefault}{\updefault}{\color[rgb]{0,0,0}$f\Big($}%
}}}}
\put(6826,-1861){\makebox(0,0)[lb]{\smash{{\SetFigFont{17}{20.4}{\familydefault}{\mddefault}{\updefault}{\color[rgb]{0,0,0}$\Big).$}%
}}}}
\put(4051,-1861){\makebox(0,0)[lb]{\smash{{\SetFigFont{17}{20.4}{\familydefault}{\mddefault}{\updefault}{\color[rgb]{0,0,0}$\Big)=f\Big($}%
}}}}
\put(5476,-1861){\makebox(0,0)[lb]{\smash{{\SetFigFont{17}{20.4}{\familydefault}{\mddefault}{\updefault}{\color[rgb]{0,0,0}$\Big)$}%
}}}}
\end{picture}%
} \end{equation*}
By a \emph{singular link} we mean a link with a finite number of self-intersections, each having distinct tangents.
A \emph{Vassiliev invariant} \cite{Vas} of type $m$ is a framed link invariant whose extension vanishes on any framed singular link with more than $m$ double points.  Let $\V_{m}$ be the vector space of all $\C$-valued Vassiliev invariants of type $m$.  The vector space $\V$ of all $\C$-valued Vassiliev invariants is filtered, 
$$\V_{0}\subset \V_{1}\subset\ldots\subset \V_{m}\subset\ldots\subset \V,$$
where $\V_{m}$ is the vector space of Vassiliev invariants of type $m$. 

Let $O^{\otimes n}$ be a disjoint union of $n$ oriented circles. 
A \emph{chord diagram} on $O^{\otimes n}$ of degree $m$ is a distinguished set of $m$ unordered pairs of points (chords), regarded up to homeomorphism preserving each connected component and the orientation.  
   Let $\chordm(O^{\otimes n})$ be the collection of chord diagram on $O^{\otimes n}$ of degree $m$, $\chordm=\bigcup_{n} \chordm(O^{\otimes n})$ and $\chord=\bigcup_{m}\chordm$. 

A \emph{weight system} of degree $m$ is a function $W: \chordm \rightarrow \C$ satisfying the four-term relation \cite[Definition 1.5]{BN}.
Let $\W_{m}$ be the vector space of all weight systems of degree $m$.  Let $\W$ be the set of all weight systems, then $\W=\oplus_{m} \W_{m}$.  

Let $\Atangle(O^{\otimes n})$ be the vector space with basis of all chord diagrams $O^{\otimes n}$ modulo the four-term relation \cite[Definition 1.7]{BN}.  Let $\Atangle_{m}(O^{\otimes n})$ be the vector space of chord diagrams of degree $m$ on $O^{\otimes n}$, modulo the four-term relation.  We also let $\Atangle=\bigcup_{n}\Atangle(O^{\otimes n})$ and $\Atangle_{m}=\bigcup_{n}\Atangle_{m}(O^{\otimes n})$.  Then $\Atangle$ is a graded vector space $\Atangle=\oplus \Atangle_{m}$.  
The space $\W_{m}$ can be represented as the linear dual of $\Atangle_{m}$.  


\subsection{Kontsevich-Bar-Natan's theorem}\label{SS:KontT}
In this subsection we state the fundamental theorem of complex-valued Vassiliev invariants.  

Let $L :\sqcup_{n}S^{1}\rightarrow\R^{3}$ be a framed singular link.  The chord diagram on $O^{\otimes n}$ of the singular link $L$ is the disjoint union of $n$ oriented circles $\sqcup_{n} S^{1}$ such that the preimages of every double point are connected by a chord.  For any Vassiliev invariant $f$ of type $m$ one can define a linear functional $[f]$ on $\chordm$ by 
$$D \mapsto f(L_{D})$$
where $L_{D}$ is any framed singular link with even framing, whose chord diagram is $D$.  Since $f$ is an invariant of type $m$ it follows that $[f]$ is well defined.  The function $[f]$ also satisfies the four-term relation allowing us to regard $[f]$ as a weight system.  Thus, we have a linear map $\V_{m}\rightarrow \W_{m}$ ($f \mapsto [f]$), whose kernel is $\V_{m-1}$.  


\begin{theorem}\label{T:BNKon}{\rm\cite{BN,Kon}}\qua
The map $\KontFunc: \V_{m}/\V_{m-1}\rightarrow \W_{m}$ is an isomorphism.
\end{theorem}

This theorem gives a purely combinatorial description of the space $\V_{0}\oplus \V_{1}/\V_{0}\oplus \V_{2}/\V_{1}\oplus\ldots$ in terms of weight systems. 
There are various proofs of this important theorem which can be found in several places including: \cite{BN,Kas,Kon}.   

The main tool Kontsevich used to prove this theorem is the construction of a family of  $\Atangle$-valued link invariant $Z_{i} \: (i=0,1\ldots)$ where $Z_{m}$ is a $\Atangle_{m}$-valued invariant of type $m$.  The function $\W_{m}\rightarrow \V_{m}$ given by $W\mapsto W\circ Z_{m}$, composed with the projection $\V_{m}\rightarrow \V_{m}/\V_{m-1}$, is the inverse $\KontFunc$.


\section{Quantum group invariants arising from Lie superalgebras of type A-G}\label{S:qg}

In this section we will define a $\C[[h]]$-valued quantum group invariant arising from a quantization of a Lie superalgebra of type A-G.  We will do this by showing that there is a functor from the category of ribbons to a category of modules over the quantization such that the restriction to links is the desired invariant.  We will also use this functor to define a (1,1)-tangle invariant.


\subsection{Ribbon categories}
We describe the concept of a ribbon category (for more details see \cite[XIV]{Kas}).  A \emph{tensor category} $\cat$ is a category equipped with a covariant bifunctor $\otimes :\cat \times \cat\rightarrow \cat$ called the tensor product, a unit object $\unit$, an associativity constraint $a$, and left and right unit constraints such that the Triangle and Pentagon Axioms are satisfied (see \cite[XI.2]{Kas}).  

When the associativity constraint and the left and right unit constraints are all identities we say the category $\cat$ is a \emph{strict} tensor category.  By Mac Lane's coherence theorem any tensor category $\cat$ is equivalent to a strict tensor category $\cat^{str}$ (see \cite[XI.5]{Kas}).  The objects of $\cat^{str}$ are finite sequences of objects of $\cat$ including the empty set.  Let $v=(V_{1}, V_{2},\dots, V_{k})$ and $v'=(V'_{1}, V'_{2},\dots, V'_{k'})$ be such sequences.  A morphism from $v$ to $v'$ in $\cat^{str}$ is a morphism from $((\cdots (V_{1}\otimes V_{2})\otimes \cdots)\otimes V_{k-1}) \otimes V_{k}$ to $((\cdots (V'_{1}\otimes V'_{2})\otimes \cdots)\otimes V'_{k'-1}) \otimes V'_{k'}$ in $\cat$ where all parentheses are placed on the left-hand side of $V_{1}$ and $V'_{1}$.  The tensor product in  $\cat^{str}$ is given by concatenation.

A \emph{braiding} on a tensor category $\cat$ consists of a family of isomorphisms $c_{V,W}: V \otimes W \rightarrow W\otimes V$, defined for each pair $(V,W)$ which satisfy the Hexagon Axiom \cite[XIII.1 (1.3-1.4)]{Kas} as well as the commutative diagram \cite[(XIII.1.2)]{Kas}.
We say a tensor category is \emph{braided} if it has a braiding.

A tensor category $\cat$ has \emph{duality} if for each object $V$ in $\cat$ there exits an object $V^{*}$ and morphisms 
$$b_{V}: \unit \rightarrow V\otimes V^{*} \;\: \text{ and } \;\: d_{V}: V^{*}\otimes V \rightarrow \unit$$
satisfying relation \cite[(XIV.2.1)]{Kas}.  A \emph{twist} in a braided tensor category $\cat$ with duality is a family $\theta_{V}:V\rightarrow V$ of natural isomorphism, defined for each object $V$ of $\cat$ satisfying relations \cite[(XIV.3.1-3.2)]{Kas}.

\begin{Df}
A (strict) ribbon category is a (strict) braided tensor category with duality and a twist. 
\end{Df}

\subsection{The category of ribbons $\rib$} We now define the category of ribbons $\rib$.  The objects of $\rib$ are finite sequences $\epsilon=(\epsilon_{1},\epsilon_{2},\ldots,\epsilon_{n})$ of $\pm$ signs.  The morphism of $\rib$ are isotopy classes of framed tangles.  A morphism from $\epsilon$ to $\epsilon'$ is an isotopy class of a framed tangle such that for any representitive $T$ of this class we have $s(T)=\epsilon$ and $b(T)=\epsilon'$ where $s$ and $b$ are the sequences defined in \S\ref{SS:klTangles}.  The composition of two morphisms $T$ and $T'$ is obtained by placing $T'$ on top of $T$.  The tensor product on objects is given by concatenation of sequences.  The tensor product of two morphisms $L$ and $L'$ is the isotopy class of the oriented framed tangle
obtained by placing $T'$ to the right of $T$.  The category $\rib$ is a strict ribbon category where braiding $c_{\epsilon,\epsilon'}$, duality $b_{\epsilon}$ and $d_{\epsilon}$, and twist $\twist_{\epsilon}$ are represented by the following framed tangles, with the blackboard framing,

\vspace{10pt}
\centerline {
\includegraphics[width=7cm]{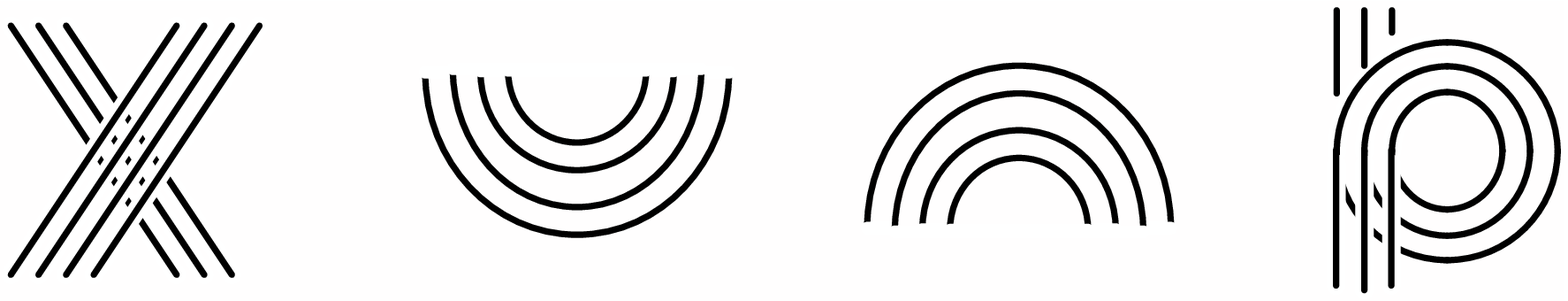}}
\noindent
respectively.
For example, the 
$$\includegraphics[width=2cm]{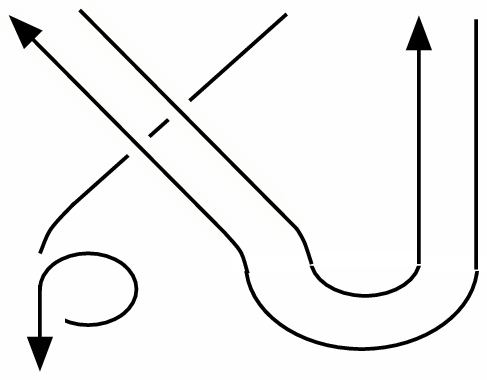}$$
is the ribbon $(c^{-1}_{(-,+),(+)}\otimes Id_{(-)}\otimes Id_{(+)})(\twist_{(+)}\otimes b_{(-,+)})$ with the blackboard framing.

We end this subsection with the following important theorem whose proof can be found in \cite{RT} (also see \cite{Kas}).

\begin{theorem}\label{T:RibFun}
Let $\cat$ be a strict ribbon category and $V$ be an object of $\cat$.  Then there exists a unique functor $F_{V}:\rib \rightarrow \cat$ which preserves the tensor product, the braiding, the duality and the twist, such that $F_{V}(+)=V$. 
\end{theorem}

\subsection{The quantum group invariant $Q_{\g,V}$}\label{SS:QGI}
Let $\g$ be Lie superalgebra of type A-G.  In this subsection we will give the definition of the quantum group invariant associated to every finite dimensional $\g$-module.  The main idea is to associate $\g$ with a Hopf superalgebra whose category of modules is a ribbon category and then apply Theorem \ref{T:RibFun}.

Let $\Uh$ be the Drinfeld-Jimbo type quantization $\g$ given by Yamane \cite{Yam94}.  It should be noted that different Dynkin diagram give different relations for $\Uh$.  For this reason we use the standard Dynkin diagram (see Table VI of \cite{K}) to construct $\Uh$.  One can also see \cite{G1} for a definition of $\Uh$.  

As a vector space $\Uh$ is isomorphic to $U(\g)[[h]]$.  Moreover, $(\Uh,R)$ is a braided Hopf superalgebra, i.e.\ a Hopf superalgebra with a universal $R$-matrix $R$.  Let  $U_{h}(\g)\text{-}Mod_{fr}$ be the category of topologically free $\Uh$-modules of finite rank (i.e.\ modules of the form $V[[h]]$ where $V$ is a $U(\g)$-module).  This category is braided where the  braiding $c_{V,W}$ for $V,W\in U_{h}(\g)\text{-}Mod_{fr}$ is defined to be
$$c_{V,W}(v\otimes w)=\flip_{V,W}\left(R(v\otimes w)\right)$$
where $v\in V$, $w\in W$ and $\flip$ is given in (\ref{E:Flip}).  We want $U_{h}(\g)\text{-}Mod_{fr}$ to be a ribbon category.  The duality is defined by the evaluation and coevaluation morphisms.  Next we define an element in $\Uh$ which we use to define a twist in $U_{h}(\g)\text{-}Mod_{fr}$.

Let $R=\sum \alpha_{i}\otimes \beta_{i}$ be the $R$-matrix in $\Uh$.  We define the elements $u, \theta_{2} \in \Uh$ by 
$$u= \sum (-1)^{\p \alpha_{i} \p \beta_{i}}S(\beta_{i})\alpha_{i} \: \: \text{ and } \:\: \theta_{2}=S(u)u.$$ 
Since $R\equiv 1 \otimes 1 \mod h$ we have that $u \equiv \theta_{2 } \equiv 1  \mod h$.   Therefore, $\theta_{2}$ has a square root which we denote by $\theta$.  For $V\in U_{h}(\g)\text{-}Mod_{fr}$ define $\twist_{V}: V\rightarrow V$ by
$$\twist_{V}(v)=\theta^{-1}v$$
for $v\in V$.  The family  $(\twist_{V})_{V}$ is a twist in $U_{h}(\g)\text{-}Mod_{fr}$ (for details see \S 4.1 of \cite{Oht:QInv}).  

Therefore, $U_{h}(\g)\text{-}Mod_{fr}$ is a ribbon category and so for each finite dimensional $\g$-module $V$, Theorem \ref{T:RibFun} implies the existence of a functor 
$$F_{\widetilde{V}}:\rib \rightarrow (\Uh\text{-}Mod_{fr})^{str}$$
 such that $F_{\widetilde{V}}(+)=(\widetilde{V})$ where $\widetilde{V}:=V[[h]]$.
The endomorphism ring of the identity element in $(\Uh\text{-}Mod_{fr})^{str}$ is $(\C[[h]])$ which is isomorphic to $\C[[h]]$.  Using this isomorphism, we define the $\C[[h]]$-valued framed tangle invariant $Q_{\g,V}$ to be the restriction of $F_{\widetilde{V}}$ to tangles.  We call $Q_{\g,V}$ the $\C[[h]]$-valued Reshetikhin-Turaev quantum group invariant of framed tangles arising from $\Uh$ and $V$.

When $\g$ is a semisimple Lie algebra and $\Uh$ is the D-J algebra, Lin \cite{Lin} showed that the $m$th coefficient of $Q_{\g,V}$ is a Vassiliev invariant of type $m$.  Lin's techniques also show that  the $m$th coefficient of $Q_{\g,V}$ is an invariant of type $m$, in the case when $\g$ is a Lie superalgebra of type A-G.  This fact is also a direct consequence of Theorem \ref{T:BNKon} and Theorem \ref{T:KontInv} which we prove in \S\ref{SS:KontInv}.  We summarize in this subsection with the following theorem.

\begin{theorem}\label{T:qg}
Let $\g$ be Lie superalgebra of type A-G and let $V$ be a finite dimensional $\g$-module.  Then there exists a $\C[[h]]$-valued R-T quantum group invariant $Q_{\g,V}$ such that the $m$th coefficient of $Q_{\g,V}$ is a Vassiliev invariant of type $m$. 
\end{theorem}               
                                                         
\subsection{The (1,1)-tangle quantum group invariant $\T{Q}_{\g, \Vhat}$}\label{SS:TQGI}
When the dimensions of the even and odd subspace of a $\g$-module $V$ are equal the invariant $Q_{\g,V}$ is zero.  In such a situation the standard technique is to consider (1,1)-tangles.  In which case then the Reshetikhin-Turaev quantum group construction associates each (1,1)-tangle with an endomorphism of $\widetilde{V}$.  When $V$ is simple weight module or more generally when the $\g$-invariant endomorphisms of $V$ are one dimensional the R-T construction can be interpreted as a $\C[[h]]$-valued (1,1)-tangle invariant.  In this subsection we will explain this construction.

First, we need a technical observation.  Let $V[[h]]$ and $W[[h]]$ be topologically free $U(\g)[[h]]$-modules.  It is well known that any morphism $f:V[[h]]\rightarrow W[[h]]$ is determined by $f|_{V}:V\rightarrow W[[h]]$.  This implies that 
$$\End_{U(\g)[[h]]}(V[[h]])\cong \End_{U(\g)}(V)[[h]].$$
In the rest of this paper we will use this isomorphism to identify these spaces.

 Let $\Vhat$ be a finite dimensional module of a Lie superalgebra $\g$ of type A-G such that $\g$-invariant endomorphisms of $\Vhat$ are one dimensional, i.e.\ $\End_{U(\g)}(\Vhat)^{\g}\cong \C$.   Note that if $\Vhat$ is a highest weight module by Schur's Lemma we have $\End_{U(\g)}(\Vhat)^{\g}\cong \C$. 
Let $T$ be a framed $(1,1)$-tangle.  Since all the building blocks of $F_{\Vhat[[h]]}$ are $\g$-invariant we have 
$$F_{\Vhat[[h]]}(T)\in (\End_{U(\g)}(\Vhat)^{\g})[[h]].$$
Since $\End_{U(\g)}(\Vhat)^{\g}\cong \C$ we can identify $F_{\Vhat[[h]]}(T)$ with a scalar in $\C[[h]]$.  We define $\T{Q}_{\g, \Vhat}(T)$ to be this power series.


\section{Weight systems arising from Lie superalgebras}\label{S:ws}

Let $\g$ be finite dimensional Lie superalgebra with a non-degenerate supersymmetric invariant even 2-tensor $t \in \g \otimes \g$.  
Let $V$ and $\Vhat$ be finite dimensional $\g$-modules such that $\End(\Vhat)^{\g}\cong \C$.   
In this section, we will define a family of weight systems $W_{\g,V}$ associated to the pair $(\g,V)$.  We will also define a family of weight systems $\T{W}_{\g,\Vhat}$ on $(1,1)$-tangles (see Definition \ref{D:chorddiagTangle}).  

We construct $W_{\g, V}$ in two steps: 1) show that the category of finite dimensional $\g$-modules $U(\g)\text{-}Mod_{f}$ is a infinitesimal symmetric category, 2) show that the category of chord diagrams on tangles $\A$ is the universal infinitesimal symmetric category and so there exists a functor $\A\rightarrow U(\g)\text{-}Mod_{f}$ which gives the desired weight system when restricted to chord diagrams on links.

\subsection{Infinitesimal symmetric categories}
A braided tensor category $\sym$ is \emph{symmetric} if the braiding $(\sigma_{V,W})_{V,W}$ satisfies 
$$\sigma_{V,W}\circ \sigma_{W,V}=id_{V\otimes W}$$
 for all pairs $(V,W)$ of object in $\sym$. 
Let $\sym$ be a strict symmetric $\C$-linear tensor category.  An \emph{infinitesimal braiding} on $\sym$ is a family $(t_{V,W})_{V,W}$ of functorial endomorphisms of $V\otimes W$, defined for all pairs $(V,W)$ of objects of $\sym$ such that 
\def\strutt{\vrule width0pt depth 10pt height 8pt}
\begin{gather}
\label{E:SymInf1}
\strutt\sigma_{V,W} \circ t_{V,W}= t_{W,V}\circ \sigma_{V,W},\\ 
\label{E:SymInf2}
\strutt t_{U,V\otimes W}=t_{U,V}\otimes id_{W} + (\sigma_{U,V}\otimes id_{W})^{-1}\circ (id_{V}\otimes t_{U,W})\circ (\sigma_{U,V}\otimes id_{W})
\end{gather}
for all objects $U,V,W$ is $\sym$.
 A category $\sym$ as above is an \emph{infinitesimal symmetric category} if it has an infinitesimal braiding.  

\begin{prop}\label{P:UgInf}
The category $U(\g)\text{-}Mod_{f}$ of finite dimensional $\g$-modules is a infinitesimal symmetric category where the symmetric braiding and infinitesimal braiding are given by 
\begin{gather}
\strutt\sigma_{V,W}=\flip_{V,W}\notag\\
\label{E:DefInfinBraiding}
\strutt t_{V,W}(v\otimes w)=t(v \otimes w)
\end{gather}
for $V,W\in U(\g)\text{-}Mod_{f}$, $v \in V$ and $w \in W$.  Recall the $\flip$ is defined in \S\ref{SS:SLS}.
\end{prop}
\begin{proof}
 The invariance of $t$ implies that $t_{V,W}$ is an endomorphism of $V\otimes W$.  A straight forward calculation shows that $t_{V,W}$ is functorial.   Since $t$ is super-symmetric we have $\sigma(t)=t$.  This implies relation (\ref{E:SymInf1}).  Finally, let us show that relation (\ref{E:SymInf2}) holds.    
  Let $t=\sum m_{i}\otimes p_{i}$.  For $u \in U$, $v \in V$ and $w\in W$ we have 
 \begin{align*}
\label{}
  t_{U, V\otimes W}(u\otimes v\otimes w)
      &=((1\otimes \Delta)(t))(u\otimes v\otimes w)   \\
      &=(\sum m_{i}\otimes (p_{i}\otimes 1 +1\otimes p_{i}))(u\otimes v\otimes w)\\
      & = \sum \left( (-1)^{\p u \p{p}_{i}}m_{i}u\otimes p_{i}v\otimes w +(-1)^{\p u \p{p}_{i} + \p v  \p{p}_{i}} m_{i}u\otimes v \otimes p_{i}w \right)
\end{align*} 
 and 
  \begin{align*}
\label{}
  \big(& t_{U,V}\otimes  id_{W} + (\sigma_{U,V} \otimes id_{W})^{-1}\circ (id_{V}\otimes t_{U,W})\circ (\sigma_{U,V}\otimes id_{W}) \big) (u\otimes v\otimes w) \\
   &= (t_{U,V}(u\otimes v))\otimes w + (\sigma_{U,V}\otimes id_{W})^{-1}\circ (id_{V}\otimes t_{U,W})((-1)^{\p u \p v}v\otimes u \otimes w)  \\
    &  =\sum (-1)^{\p u \p{p}_{i}}m_{i}u \otimes p_{i} v \otimes w + (\sigma_{U,V}\otimes id_{W})^{-1}\left( \sum (-1)^{\p u \p v + \p{p}_{i}\p u} v\otimes m_{i}u \otimes p_{i}w \right) \\
    & =  \sum \left( (-1)^{\p u \p{p}_{i}}m_{i}u \otimes p_{i} v \otimes w+   (-1)^{\p u \p v + \p{p}_{i}\p u + \p v \p{m}_{i} + \p v \p u} m_{i}u\otimes v \otimes p_{i}w\right)
\end{align*}
where $\p u \p v + \p v \p u = \p 0$ and $\p{m}_{i}=\p{p}_{i}$ (since $t$ is even).  Therefore, relation (\ref{E:SymInf2}) follows.
\end{proof}
\subsection{The category $\A$}

\begin{Df}\label{D:chorddiagTangle}
Let $T$ be a framed tangle.
A chord diagram on $T$ of degree $m$ is the tangle $T$ with a distinguished set of $m$ unordered pairs of points of $T \backslash \partial T$, considered up to homeomorphisms preserving each connected component and the orientation.  
\end{Df}

Let $\Etangle(T)$ be the vector space with basis given by all chord diagrams on $T$.  If $f:T\rightarrow T'$ is a homeomorphism of tangles, then $f$ induces a isomorphism $\Etangle(T)\cong \Etangle(T')$.  Since any tangle is isomorphic to the ``unknotted'' tangle, the isomorphism class of $\Etangle(T)$ only depends on the number of line segments and circles which make up $T$.

Let $\Atangle(T)$ be the vector space defined as the quotient of $\Etangle(T)$ by the four term relation.  The vector space $\Atangle(T)$ has a natural grading, $\Atangle(T)=\oplus_{m\geq 0} \Atangle_{m}(T)$,
where $\Atangle_{m}(T)$ is all chord diagram on $T$ of degree $m$. 

Next we describe the category of chord diagrams $\A$.  As in the category of ribbons $\rib$ the objects of $\A$ are finite sequences of $+$ and $-$, along with $\emptyset$.  The morphisms of $\A$ are elements of $\Atangle(T)$ for some framed tangle $T$.  The tensor structure is defined as in the category $\rib$.   Moreover, the braiding of the category $\rib$ induces a braiding in $\A$.  Since chord diagrams are considered up to homeomorphism this braiding is symmetric, i.e.
$${\includegraphics[width=4cm]{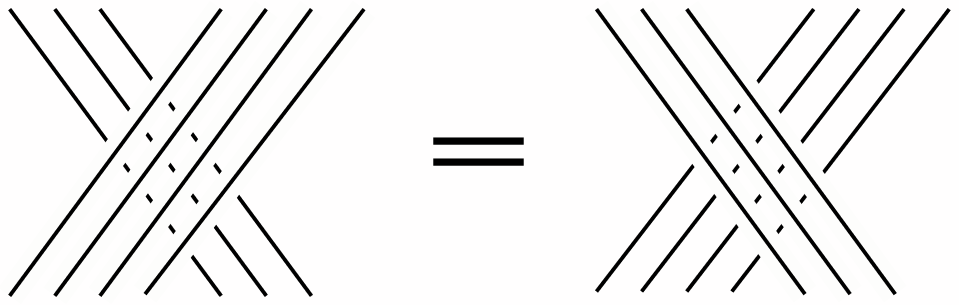}}$$
We will now describe the infinitesimal braiding on $\A$.  Let $\epsilon$ and $\epsilon'$ be objects of $\A$, i.e.\ $\epsilon=(\epsilon_{1},\ldots,\epsilon_{n}) $ where $\epsilon_{i}=\pm$.  Let $Id_{\epsilon}$ be the identity of $\epsilon$.  In other words, $Id_{\epsilon}$ is the oriented ($n$,$n$)-tangle given by $n$ vertical (parallel) strings whose orientations are given by $\epsilon$ (see \S\ref{SS:klTangles}).  Define $t_{\epsilon,\epsilon'}$ to be the sum of all chord diagrams with exactly one chord on $Id_{\epsilon}\otimes Id_{\epsilon'}$ such that the one of the end points of the chord is on $Id_{\epsilon}$ and the other is on $Id_{\epsilon'}$.  

The following theorem can be found in \cite[XX]{Kas}.
\begin{theorem}\label{T:Astrict}
With the above symmetric braiding and infinitesimal braiding the category $\A$ is a strict infinitesimal symmetric category with duality.
\end{theorem}

The category $\A$ is the universal infinitesimal symmetric category.  This statement is made precise by the following theorem whose proof can be found in \cite[XX.8]{Kas}.

\begin{theorem}\label{T:FunctorAtoSym}
Let $\sym$ be an infinitesimal symmetric category with duality and let $V$ be an object of $\sym$.  Then there exists a unique functor $F_{V}:\A \rightarrow \sym$ preserving the tensor product, symmetry, infinitesimal braiding and the duality such that $F_{V}(+)=V$.
\end{theorem}

Recall that from Proposition \ref{P:UgInf} we have that $U(\g)\text{-}Mod_{f}$ is a infinitesimal symmetric category.  Therefore, for each finite dimensional $\g$-module $V$, Theorem \ref{T:FunctorAtoSym} implies there exists a functor
\begin{equation}
\label{E:FuncG}
G_{\g,V}:\A \rightarrow U(\g)\text{-}Mod_{f}
\end{equation}
which preserving the tensor product, symmetry, infinitesimal braiding and the duality such that $G_{\g,V}(+)=V$.

Next we will make two observation and then we will give the definition of the weight system $W_{\g,V}$.  First, if $D$ is a chord diagram on a link then $ G_{\g,V}(D)$ is an endomorphism of the identity $I=\C[[h]]$, where $\End(I)\cong \C[[h]]$.  We use this isomorphism to identify $ G_{\g,V}(D)$ with a scalar.  
Second, let $\Atangle(\downarrow)$ be the vector space of chord diagrams on (1,1)-tangles modulo the four term relation.  Let $D \in \Atangle(\downarrow)$ then $G_{\g,V}(D)$ is a endomorphism of $V$.  Moreover, since all of the morphims used to define a infinitesimal symmetric category are functorial we have that $G_{\g,V}(D)$ is a $\g$-invariant endomorphism.  

\begin{Df}\label{D:ws}
Let $\g$ be finite dimensional Lie superalgebra with a non-degenerate supersymmetric invariant even 2-tensor $t \in \g \otimes \g$ and let $V,\Vhat$ be finite dimensional $\g$-modules such that $\End(\Vhat)^{\g}\cong \C$.  Define $W_{\g,V}$ to be the linear functional on $\Atangle$ given by the restriction of $G_{\g,V}$ to $\Atangle$. (Recall $\Atangle$ is the space of chord diagrams on links modulo the four-term relation.) Define $\T{W}_{\g,\Vhat}$ to be the linear functional on $\Atangle(\downarrow)$ given by $\T{W}_{\g,\Vhat}(D)=G_{\g,\Vhat}(D)\in \End(\Vhat)^{\g}\cong \C$ for $D$ in $\Atangle(\downarrow)$. 
\end{Df}

\begin{remark}\label{R:ws}
(1)\qua In the case when $\g$ is a finite dimensional Lie algebra the family of weight systems $W_{\g,V}$ becomes Bar-Natan's construction \cite{BN}.  

(2)\qua The linear functional $\T{W}_{\g,\Vhat}$ restricted to one component (1,1)-tangles becomes the family of weight systems defined by Vaintrob \cite{V}. 
\end{remark}

\section{The Kontsevich integral and invariants arising from Lie superalgebras}\label{S:gl}

In this section we will define the Kontsevich integral and then state and prove the main results of this paper.  

\subsection{The Kontsevich integral}
Here using algebraic techniques we will define an invariant whose restriction to links is the  Kontsevich integral.  For this reason we call the invariant the Kontsevich integral.

Let $\sym$ be an infinitesimal symmetric category with duality.  Let $\sym[[h]]$ be the category whose objects are the same as the objects of $\sym$ and morphisms are given by
$$\Hom_{\sym[[h]]}(V,W):= \Hom_{\sym}(V,W)[[h]].$$  By the following theorem one can extend $\sym$ to a strict ribbon category $\sym[[h]]^{str}$.

\begin{theorem}\label{T:SymIsBraided}
Let $\sym$ be a infinitesimal symmetric category with symmetry 
$(\sigma_{V,W})_{V,W}$ and infinitesimal braiding $(t_{V,W})_{V,W}$.  The tensor structure of $\sym$ extends to a braided tensor category on $\sym[[h]]$ such that the associativity constraint $a$ and the braiding $c$ are given by 
\begin{equation}
\label{EAssocBraiding}
a_{U,V,W}=\Phi(t_{U,V},t_{V,W}) \hspace{10pt} \text{and}\hspace{10pt} c_{V,W}=\sigma_{V,W} \circ e^{ht_{V,W}/2}
\end{equation}
where $U,V,W$ are any objects in $\sym$ and $\Phi$ is the Drinfeld associator.
Moreover, if the category $\sym$ has a duality then $\sym[[h]]^{str}$ is a strict ribbon category. 
\end{theorem}
\begin{proof}
See \cite{Kas} theorem XX.6.1 and theorem XX.7.1.
\end{proof}
From Theorems \ref{T:Astrict} and \ref{T:SymIsBraided} we have that $ \A[[h]]^{str}$ is a strict ribbon category.  Applying Theorem \ref{T:RibFun} to the category $\A[[h]]^{str}$ and the object $(+)$ we have a framed tangle invariant
$$Z: \rib \rightarrow \A[[h]]^{str}$$
where $Z(+)=(+)$.  We set $Z(L)=\sum_{i \geq 0}Z_{i}(L)h^{i}$.  

We will now derive an important property of $Z$.  Let $T: D\rightarrow \R^{3}$ be a singular tangle (with $m$ double points) where $D$ is a one dimensional smooth compact manifold.  The chord diagram associated to the singular tangle $T$ is $D$ such that the preimage of every double point is connected by a chord. 

When evaluating $Z$ each double point of a singular link locally contributes $c_{\epsilon,\epsilon'}- c_{\epsilon',\epsilon}^{-1}$.  From the definition of $Z$ the braiding is given by (\ref{EAssocBraiding}) and so for every singular point $Z$ has a local contribution of 
$$e^{ht_{\epsilon,\epsilon'}/2}-e^{-ht_{\epsilon,\epsilon'}/2}=ht_{\epsilon,\epsilon'} + \text{ higher order terms}.$$
From the last equality it follows that 
\begin{equation}
\label{E:ZLD}
Z(T)\equiv D h^{m} \mod h^{m+1}.
\end{equation}
In particular, when $D$ is a chord diagram of a link with $m$ chords and  $L_{D}$ is a singular link whose chord diagram is $D$ we have 
\begin{align}
\label{E:ZmLD}
   Z_{m}(L_{D})=D.   
\end{align}
From equality (\ref{E:ZmLD}) it follows that 
$$[w \circ Z_{m}](D)=(w \circ Z_{m})(L_{D})=w(D)$$
where $[\cdot]$ is the linear map defined is \S\ref{SS:KontT}.  Thus, the linear map 
 $\W_{m} \rightarrow \V_{m}$ given by $w \mapsto w \circ Z_{m}$ composed with the projection $\V_{m}\rightarrow \V_{m}/\V_{m-1}$ is inverse of $\KontFunc$, the linear map given in Theorem \ref{T:BNKon}.  For this reason we call $Z$ the Kontsevich integral.  

 \subsection{Main Theorem}\label{SS:KontInv}

\begin{theorem}\label{T:KontInv}
Let $\g$ be a Lie superalgebra of type A-G.  Let $V,\Vhat$ be finite dimensional $\g$-modules such that $\End(\Vhat)^{\g}\cong \C$.  Then
\begin{equation}
\label{E:KontInv1}
Q_{\g,V}=W_{\g,V}\circ Z : \{\text{framed links}\} \rightarrow \C[[h]]
\end{equation}
and
\begin{equation}
\label{E:KontInv2}
\T{Q}_{\g,\Vhat}(T)= \T{W}_{\g,\Vhat}\circ Z: \{\text{framed (1,1)-tangles}\} \rightarrow \C[[h]].
\end{equation}
\end{theorem}


\begin{proof}
Recall that $Q_{\g,V}$ is the restriction of $F_{\widetilde{V}}$ to framed links.  To prove the theorem we will define two functors
$$E:(U(\g)\text{-}Mod_{f})[[h]]^{str} \rightarrow U_{h}(\g)\text{-}Mod_{fr},$$
$$G_{\g,V}[[h]] : \A[[h]]^{str}\rightarrow (U(\g)\text{-}Mod_{f})[[h]]^{str},$$
such that the restriction of $E \circ G_{\g,V}[[h]] \circ Z$ to framed links is the right side of equation (\ref{E:KontInv1}).  

We will now define two functors which we will use to define the functor $E$.  Let $t\in \g\otimes \g$ be the 2-tensor associated to the unique non-degenerate supersymmetric $\g$-invariant bilinear form on $\g$ (see \S\ref{SS:SLS}).  Let $\kzt$ be the braided quasi-Hopf superalgebra constructed in \cite{G1}.  As a vector space $\kzt$ is isomorphic to $U(\g)[[h]]$. The definition of $\kzt$ is based on the theory of the Knizhnik-Zamolodchikov differential equations.  In particular, the associator of the quasi-Hopf superalgebra is Drinfeld associator.  The standard categorical argument (see \cite[Prop. XV.1.2]{Kas}) shows that $\kzt$-$Mod_{fr}$ is a ribbon category where the associativity constraint $a$ and the braiding $c$ are given by 
$$a_{U,V,W}=\Phi(t_{12},t_{23}) \hspace{10pt} \text{and}\hspace{10pt} c_{V,W}=\flip_{V,W} \circ e^{ht/2}$$
where $U,V,W$ are any objects in $\kzt$-$Mod_{fr}$, $\Phi$ is the Drinfeld associator and $\flip_{V,W}$ is the linear map given in (\ref{E:Flip}).  

In \cite{G1} the author constructs a functor $D:  \kzt\text{-}Mod_{fr} \rightarrow U_{h}(\g)\text{-}Mod_{fr}$ which is an equivalence of categories which preserves the ribbon structures such that $D(V[[h]])=V[[h]]$.  Such a functor was first constructed by Drinfeld in the case when $\g$ is a semi-simple Lie algebra.  The techniques used by Drinfeld do not work for Lie superalgeras.  The construction given in \cite{G1} is based on the Etingof-Kazhdan quantization of Lie superbialgebras.

From Proposition \ref{P:UgInf} we have $U(\g)\text{-}Mod_{f}$ is a infinitesimal symmetric category.  By Theorem \ref{T:SymIsBraided}, $(U(\g)\text{-}Mod_{f})[[h]]$ is a ribbon category where the associativity constraint and the braiding are given by (\ref{EAssocBraiding}).  By construction the tensor categories $(U(\g)\text{-}Mod_{f})[[h]]$ and $\kzt \text{-}Mod_{fr}$ have the same associativity constraint and braiding.  Therefore the functor
$$(U(\g)\text{-}Mod_{f})[[h]] \rightarrow \kzt \text{-}Mod_{fr} $$
mapping $V \mapsto V[[h]]$ is an equivalence of categories which preserves the ribbon structure.  

We define the functor $E$ by the composition of the following equivalences:
$$(U(\g)\text{-}Mod_{f})[[h]]^{str} \rightarrow (\kzt \text{-}Mod_{fr})^{str} \xrightarrow{} \kzt \text{-}Mod_{fr}\xrightarrow{D} U_{h}(\g)\text{-}Mod_{fr}.$$
The tensor functor $E$ preserves the ribbon structure since $D$ does.

Now we define $G_{\g,V}[[h]]$.  Recall the definition of the functor $G_{\g,V}$ given in \S \ref{E:FuncG}.  From Theorem \ref{T:SymIsBraided} the functor $G_{\g,V}$ extends to a functor $$G_{\g,V}[[h]] : \A[[h]]^{str}\rightarrow (U(\g)\text{-}Mod_{f})[[h]]^{str}$$ such that $G_{\g,V}[[h]](+)=V[[h]]=\widetilde{V}$ and $G_{\g,V}[[h]]$ preserves the ribbon structure.  Let $\sum D_{i}h^{i}\in \A[[h]]$.  Then by definition 
\begin{equation}
\label{E:GVhisGVh}
G_{\g,V}[[h]]( \sum_{i\geq 0}D_{i}h^{i}) = \sum_{i\geq 0}G_{\g,V}(D_{i})h^{i}
\end{equation} 
$$E \circ G_{\g,V}[[h]] \circ Z:\rib \rightarrow U_{h}(\g)\text{-}Mod_{fr}
\leqno{\hbox{The functor}} 
$$ preserves the ribbon structure maps the object $(+) $ to $\widetilde{V}:=V[[h]]$. Therefore, from the uniqueness of Theorem \ref{T:RibFun} we have 
\begin{equation}
\label{E:FunctorsEqual}
F_{\widetilde{V}}=E \circ G_{\g,V}[[h]] \circ Z.
\end{equation}

We now prove the first assertion of the theorem.  Let $L$ be a link with $n$ components.  Then $Z(L)=\sum_{i\geq 0}D_{i}h^{i}$ where $D_{i}$ is a chord diagram on $L$.  We have 
\begin{align}
\label{E:EGW}
  (E \circ  G_{\g,V}[[h]]\circ Z)(L)&= (E \circ  G_{\g,V}[[h]])( \sum_{i\geq 0}D_{i}h^{i}) \notag \\
  &= E(\sum_{i\geq 0}G_{\g,V}(D_{i})h^{i}) \notag
   \\
    &=  \sum_{i\geq 0}W_{\g,V}(D_{i})h^{i} \notag\\
    &=  (W_{\g,V}\circ Z)(L)
\end{align}
where the first and forth equalities are definitions, the second equality follows from (\ref{E:GVhisGVh}), and the third equality follows from definition of the weight system and the fact that $E$ preserves the ribbon structure.   
The first assertion of the theorem follows from (\ref{E:EGW}) and the uniqueness of equality (\ref{E:FunctorsEqual}).  

Similarly, for a (1,1)-tangle $T$ we have 
$$(E \circ  G_{\g,\Vhat}[[h]]\circ Z)(T)=(\T{W}_{\g,\Vhat}\circ Z)(T).$$ This equality combined with the uniqueness of equality (\ref{E:FunctorsEqual}) imply the second assertion of the theorem.
\end{proof}


\begin{corollary}\label{C:UnderQ}
Let $Q_{\g,V}=\sum Q_{\g,V,m}h^{m}$ and $W_{\g,V}=(W_{\g,V,m})_{m}$. 
For all $m$, we have the weight system corresponding to  $Q_{\g,V,m}$  is $W_{\g,V,m}$, i.e. 
\begin{align*}
\label{}
   [Q_{\g,V,m}]=&W_{\g,V,m}
\end{align*}
where $[\cdot]$ is the map defined in \S \ref{SS:KontT}. 
\end{corollary}

\begin{proof}
Let $D$ be a chord diagram of a link with $m$ chords and let $L_{D}$ be a framed singular link with even framing, whose chord diagram is $D$.
We have 
\begin{align*}
 [Q_{\g,V,m}](D) &=  Q_{\g,V,m}(L_{D})  \\ & = W_{\g,V,m} (Z_{m}(L_{D}))  \\
      &   = W_{\g,V,m}(D)
\end{align*}  
where the first equality is by definition, the second follows from Theorem \ref{T:KontInv} and the third follows from equity (\ref{E:ZmLD}).
\end{proof}

Let $W=(W_{m})_{m\in \Z}$ be a family of weight systems and $Q=\sum Q_{m}h^{m}$ be a quantum group invariant.  If $[Q_{m}]=W_{m}$, for all $m$, we say the weight systems corresponding to $Q$ are equal to the family $W$.
We summarize this subsection with the following theorem:

\begin{theorem}\label{T:gl}
Let $V$ be a finite dimensional module of a Lie superalgebra $\g$ of type A-G.  Let $Q_{\g,V}$ be the quantum group invariant of Theorem \ref{T:qg}, let $W_{\g,V}$ be the family of  weight system of Definition \ref{D:ws} and let $Z$ be the Kontsevich integral $Z$.  Then $Q_{\g,V}=W_{\g,V}\circ Z$.  Moreover, the weight systems corresponding to $Q_{\g,V}$ are  equal to the family $W_{\g,V}$.  In summary the following diagram commutes
\begin{equation}\label{D:gV}
\xymatrix{     & \{ \g,V\}\ar[ld]_{\text{ Definition \ref{D:ws} }} \ar[rd]^{\text{ Theorem \ref{T:qg} }} & \\                    
                     \W  \ar@<-2.5pt>[rr]_{Z^{*}}& & \V \ar@<-1.5pt>[ll]_{[\cdot]} }
\end{equation}
where $[\cdot]$ is the map defined in subsection \ref{SS:KontT}. 
\end{theorem}

\subsection{The Kontsevich integral and (1,1)-tangle invariants}
Not surprisingly the theory of Vassiliev invariants can be formulated in the case of (1,1)-tangles.  In this subsection we do this and then prove a theorem analogous to Theorem \ref{T:gl}  for (1,1)-tangles.

 A \emph{Vassiliev invariant of (1,1)-tangle of type $m$} is a framed (1,1)-tangle invariant whose extension vanishes on any framed singular (1,1)-tangle with more than $m$ double points.  Let $\T{\V}_{m}$ be the vector space of all $\C$-valued Vassiliev (1,1)-tangle invariants of type $m$.  A \emph{(1,1)-tangle weight system} of degree $m$ is a linear functional on the space $\Atangle_{m}(\downarrow)$ of all chord diagrams on $\downarrow$ modulo the four-term relation.  Let $\T{\W}_{m}$ be the vector space of all (1,1)-tangle weight system of degree $m$.  When it is clear we will omit the word (1,1)-tangle when referring to Vassiliev invariants and weight systems.

For a Vassiliev invariant $f$ of type $m$ we define $[f]$ to be the weight system given by $[f](D)=f(T_{D})$ where $T_{D}$ is a (1,1)-tangle with even framing, whose chord diagram is $D$.   
The assignment $f \mapsto [f]$ defines a linear map $[\cdot]:\T{\V}_{m}\rightarrow \T{\W}_{m}$ whose kernel is $\T{\V}_{m-1}$.  Therefore, $[\cdot]$ induces a linear map 
$$\T{\V}_{m}/\T{\V}_{m-1}\rightarrow \T{\W}_{m}.$$
As in the case of links, the inverse of this map is the linear map $\T{\W}_{m}\rightarrow \T{\V}_{m}$ given by $W\mapsto W\circ Z_{m}$, composed with the projection $\T{\V}_{m}\rightarrow \T{\V}_{m}/\T{\V}_{m-1}$.

\begin{theorem}\label{T:QWZ11} Let $\T{Q}_{\g, \Vhat}$ be the (1,1)-tangle quantum group invariant defined in \S\ref{SS:TQGI}.  Let $\T{W}_{\g,\Vhat}$ be the weight system given in Definition \ref{D:ws}.  Then we have
\begin{enumerate}
  \item the $m$th coefficient of $\T{Q}_{\g, \Vhat}$ is a Vassiliev invariant of type $m$,\label{TI:Vas}
  \item  $\T{W}_{\g,\Vhat}=\T{Q}_{\g, \Vhat}\circ Z$,\label{TI:WZQ}
  \item the weight systems corresponding to  $\T{Q}_{\g, \Vhat}$ are the equal to the family   $\T{W}_{\g,\Vhat}$.
\end{enumerate}
\end{theorem}
\begin{proof}
Assertion (\ref{TI:WZQ}) is a restatement of equality (\ref{E:KontInv2}).  Assertion (\ref{TI:Vas}) is a direct consequence of (\ref{TI:WZQ}).  The final assertion as follows.  By definition the weight system corresponding to  $\T{Q}_{\g, \Vhat,m}$ is $[\T{Q}_{\g, \Vhat,m}]$.  Let $D$ be a chord diagram of a (1,1)-tangle with $m$ chords and let $T_{D}$ be a singular (1,1)-tangle whose chord diagram is $D$.
We have 
\begin{align*}
 [\T{Q}_{\g,\Vhat,m}](D) &=  \T{Q}_{\g,\Vhat,m}(T_{D})   = \T{W}_{\g,\Vhat,m} (Z_{m}(T_{D}))    = \T{W}_{\g,\Vhat,m}(D)
\end{align*}  
where the first equality is by definition, the second follows from Theorem \ref{T:KontInv} and the third follows from equivalence (\ref{E:ZLD}).
\end{proof}

\subsection{General theorem}\label{SS:GenThem}

In this subsection we list the data need to prove a theorem analogous to Theorem \ref{T:gl} for a general Lie superalgebras.

Let $\g$ be finite dimensional Lie superalgebra with an even element $r \in \g \otimes \g$ such that $t:=r+ \flip_{\g,\g}(r)$ is a non-degenerate supersymmetric $\g$-invariant element.  Let $\kzt$ be the braided quasi-Hopf superalgebra constructed in \cite{G1}. 
Let  $\Uh$ be a superalgebra with the following properties:
\begin{enumerate}
  \item  as a vector space $\Uh$ is isomorphic to $U(\g)[[h]]$,\label{I:UhIso}
  \item  $\Uh$ is a topological Hopf superalgebra,
  \item  $\Uh$ has a universal $R$-matrix, i.e.\ a  homogeneous element  $R\in \Uh^{\otimes 2}$ satisfying relations (2.1), (2.3) and (2.4) of \cite[VIII.2]{Kas},\label{I:UhRmatrix}
  \item the $R$-matrix is of the form $R\equiv 1\otimes 1 +rh \mod h^{2}$,\label{I:UhFormR}
  \item  $\Uh$ has a invertible homogeneous element $\theta$ satisfying relations \cite[(XIV.6.5)]{Kas}, \label{I:UhCentral}
  \item  for each finite dimensional module $V$ of $\g$ we can associate a topologically free $\Uh$-module $\widetilde{V}\cong V[[h]]$,\label{I:Uhfinitedim}
  \item there exists a functor $D:  \kzt\text{-}Mod_{fr} \rightarrow U_{h}(\g)\text{-}Mod_{fr}$ which is an equivalence of categories which preserves the ribbon structures such that $D(\widetilde{V})=\widetilde{V}$ for all $\widetilde{V} \in \kzt\text{-}Mod_{fr}$.\label{I:functorD}
\end{enumerate}

Let $V$ and $\Vhat$ be finite dimensional $\g$-modules such that $\End(\Vhat)^{\g}\cong \C$.   Let $Q_{\g,V}$ ($\T{Q}_{\g,\Vhat}$) be the $\C[[h]]$-valued Reshetikhin-Turaev quantum group invariant of framed tangles arising from $\Uh$ and $V$ (resp. $\Vhat$).  Let $W_{\g,V}$ and  $\T{W}_{\g,\Vhat}$ be the families of weight systems defined in \S \ref{S:ws}.  The following theorem follows exactly as in the proofs of Theorem \ref{T:KontInv}, Corollary \ref{C:UnderQ} and Theorem \ref{T:QWZ11}.

\begin{theorem}\label{T:GenQWZ}
Let $\g$ be a finite dimensional Lie superalgebra with an even element $r \in \g \otimes \g$ such that $t:=r+ \flip_{\g,\g}(r)$ is a non-degenerate supersymmetric $\g$-invariant element.  Suppose that there exists a superalgebra $\Uh$ and a functor $D$ satisfying assumptions (\ref{I:UhIso})-(\ref{I:functorD}).  Then the $m$th coefficients of $Q_{\g,V}$ and  $\T{Q}_{\g,\Vhat}$ are Vassiliev invariant of type $m$ such that  
\begin{align}
\label{}
   Q_{\g,V}=& W_{\g,V}\circ Z, & \T{Q}_{\g,\Vhat}=& \T{W}_{\g,\Vhat}\circ Z.
\end{align} 
Moreover, the weight systems corresponding to  $Q_{\g,V}$ ($\T{Q}_{\g,\Vhat}$) are equal to the family $W_{\g,V}$ (resp. $\T{W}_{\g,\Vhat}$).
\end{theorem}

\section{The Links-Gould invariant}\label{S:ap}

In this section we will use the results of section \ref{S:gl} to investigate the Links-Gould invariant.  First we show that invariants arising from the general linear Lie superalgebra $\gl$ can be viewed via weight systems and the Kontsevich integral.

\subsection{Invariants arising from $\gl$}
In this subsection we show that for the general linear Lie superalgebra $\gl$ there exists a superalgebra $U_{h}(\gl)$ and functor $D$ that satisfy assumptions (\ref{I:UhIso})-(\ref{I:functorD}).  Then we will conclude from Theorem \ref{T:GenQWZ} that $Q_{\gl,V}=W_{\gl,V}\circ Z$.

The general linear Lie superalgebra $\gl$ is a not a Lie superalgebra of type A-G.  However, for $m\neq n$, $\gl$ is a one-dimensional central extension of $\slmn$ which is the Lie superalgebra of type A$(m-1,n-1)$.  Now $\mathfrak{sl}(n|n)$ has a one-dimensional ideal consisting of the scalar matrices and is not the Lie superalgebra of type A$(n-1,n-1)$.  For this reason in this section we will assume $n\neq m$.

We now show that $Q_{\gl,V}$ and $W_{\gl,V}$ exist and correspond via the Kontsevich integral.  Let $t$ be the 2-tensor associated to the unique non-degenerate supersymmetric invariant bilinear form on $\slmn$.  Let $s\in \gl \otimes \gl$ be an extension of $t$ of the form 
\begin{equation}
\label{E:Defs}
s= a (I\otimes I) + t
\end{equation}
where $a$ is a non-zero complex number and $I$ is the identity element of $\gl$.  Using $s$ one can still construct a Drinfeld-Jimbo type superalgebra $U_{h}(\gl)$ \cite{GThesis}.  This superalgebra is almost identical to $U_{h}(\slmn)$ with the following two modifications: (1)  $U_{h}(\gl)$ has an additional generator $E_{0}$ such that $E_{0}$ is central and the classical limit of $E_{0}$ is the identity element of $\gl$, (2) the $R$-matrix of $U_{h}(\gl)$ is slightly modified to account for $E_{0}$ and the complex number $a$.    

Let $A_{\gl,s}$ be the braided quasi-Hopf superalgebra constructed in \cite{G1}.  As shown in  \cite{GThesis}, there is a braided tensor functor 
$$D: A_{\gl,s} \text{-}Mod_{fr} \rightarrow U_{h}(\gl)\text{-}Mod_{fr}.$$ 
After straight forward adjustments accounting for the central element $E_{0}$, the construction of this functor is almost identical to the construction for $\slmn$ given in \cite{G1}.  

Let $V$ and $\Vhat$ be finite dimensional $\gl$-modules such that $\End(\Vhat)^{\gl}\cong \C$.  Using $U_{h}(\gl)$ we can define the $\C[[h]]$-valued Reshetikhin-Turaev quantum group invariant $Q_{\gl,V}$ and $\T{Q}_{\gl,\T V}$.  Let $W_{\gl,V}$ and $\T{W}_{\gl,\Vhat}$ be the weight systems given in Definition \ref{D:ws} where we use the 2-tensor $s$.  The following theorem is a special case of Theorem \ref{T:GenQWZ}.

\begin{theorem}\label{T:QWZglmn}
The $m$th coefficients of $Q_{\gl,V}$ and $\T{Q}_{\gl,\T V}$ are Vassiliev invariants of type $m$ such that $Q_{\gl,V}=W_{\gl,V}\circ Z$ and  $\T{Q}_{\gl,\T V}=\T{W}_{\gl,\Vhat}\circ Z$.  Moreover, the weight systems corresponding to $Q_{\gl,V}$ ($\T{Q}_{\gl,\T V}$) are equal to $W_{\gl,V}$ (resp. $\T{W}_{\gl,\Vhat}$).
\end{theorem}

\subsection{The Links-Gould invariant as a power series in $h$}\label{SS:LGpower} In this subsection we will give a reinterpretation of the Links-Gould invariant with $q=e^{h/2}$.  First we will show that the invariant $Q_{\gto,\Valpha}$ is zero.  Throughout this subsection we set $\g=\gto$.

Recall that the Links-Gould invariant \cite{LG} is a quantum group invariant of unframed (1,1)-tangles arising from $U_{q}(\g)$ and a family of $2|2$-dimensional simple modules $\Valpha$, $(\alpha \in \C)$.  Thinking of $\alpha$ as a variable we have this invariant is a polynomial in two variables $q$ and $q^{\alpha}$.  

For each $\alpha\in \C$ the module $\Valpha$ has a two dimensional even and two dimensional odd decomposition.   This implies that the supertrace of any scalar endomorphism of $\Valpha$ is zero.   Let $T$ be a framed $(1,1)$-tangle and let $\overline{T}$ be the link that is the closure of $T$.  Recall from \S \ref{SS:TQGI} that $F_{\Valphahat}(T)$ is an element of $\End_{\g}(\Valpha)^{\g}[[h]]$.  By Schur's Lemma $F_{\Valphahat}(T)$ is a scalar matrix and so has zero supertrace.  Taking the closure of a $(1,1)$-tangle under $F_{\Valphahat}$ corresponds to taking the supertrace of the endomorphism associated to the tangle, in other words $F_{\Valphahat}(\overline{T})$ is the supertrace of $F_{\Valphahat}(T)$.  Thus, $Q_{\g,\Valpha}=F_{\Valphahat}|_{Links}=0$.   For this reason the Links-Gould invariant is defined on framed $(1,1)$-tangles.  

Let $s$ be the 2-tensor defined in equation (\ref{E:Defs}) with $a=2+2\alpha$.  Define $LG$ to be the (1,1)-tangle invariant $\T{Q}_{\g,\Valpha}$ arising from $s$, $U_{h}(\g)$ and $V$.  Let $\sigma$ be the matrix that represents the action of the $R$-matrix on $\Valpha[[h]]\otimes \Valpha[[h]]$.  The entries of $\sigma$ take values in $\C[\alpha][[h]]$ (see \cite[\S 4.1]{DWLK}).  Therefore, $LG$ is an $\C[\alpha][[h]]$-valued invariant.  $LG$ is the Links-Gould invariant with $q=e^{h/2}$.
Define the weight system $\lgws$ to be $\T{W}_{\g,\Vhat}$ arising from $s$, $\g$ and $\Vhat$.   

The following theorem is a special case of Theorem \ref{T:QWZglmn}.

\begin{theorem}\label{P:LGandWS}
The Links-Gould invariant $LG$ is equal to the weight system $\lgws$ composed with the Kontsevich integral $Z$.  In other words, 
\begin{align}
\label{E:LGandWS}
   LG=\lgws\circ Z :\{ \text{framed (1,1)-tangles}\} \rightarrow \C[[h]].
\end{align}
Conversely, the weight systems corresponding to $LG$ are equal to the family $\lgws$. 
\end{theorem}  

We will now show that $LG$ defines an invariant of knots.  It should be noted this fact also follows from the observation that isotopy classes of knots are in one to one correspondence with isotopy classes of long knots \cite{TUR}.  Let $B(\downarrow)$ ($B(O)$) be the space of chord diagrams on the ``unknoted'' one component (1,1)-tangle (resp. on the circle) modulo the four term relation.  It is well known that the map $B(O)\rightarrow B(\downarrow)$ given by slitting the circle at a point (which is not an endpoint of a chord) is a well defined isomorphism (see \cite{BN,Kas}).  In other words, if $D$ and $D'$ are two chord diagram coming from slitting a chord diagram on the circle, then $D=D'$ in $B(\downarrow)$.

Let $K$ be a knot.  By slitting $K$ at some point $p$ we get a (1,1)-tangle $T$.  Let $T'$ be another (1,1)-tangle gotten from slitting $K$ at some point, possibly different than $p$.  From the discussion above (or more precisely \cite[Lemma XX.3.1]{Kas}) we have $Z(T)=Z(T')$.  Thus, we are led to the following definition.

\begin{Df}
Let $K$ be a knot and $T$ a (1,1)-tangle created from slitting $K$ at some point.  Define
$$LG(K)=\lgws\circ Z(T).$$
\end{Df}

From Theorem \ref{P:LGandWS} we have that $LG(T)=\lgws\circ Z(T)$.  Therefore, the definition of $LG(K)$ is natural.

\begin{remark}\label{R:FramingIndepOfLG}
(1)\qua From Theorem \ref{P:LGandWS} it follows that for each (1,1)-tangle chord diagram $D$ we have $\lgws(D)$ is a polynomial in $\alpha$.  

(2)\qua As stated above $LG$ is a (1,1)-tangle invariant, so it is framing independent.  By direct calculation one can show that $\lgws$ satisfies the one term relation meaning that it is zero on diagrams with an isolated chord.  It follows that the invariant $ \lgws\circ Z$ is framing independent.   Thus, equality (\ref{E:LGandWS}) is framing independent.  Note the framing independence of both sides of equality (\ref{E:LGandWS}) follows from the choice of the scalar $a$ in $s$.
\end{remark}

\subsection{Cablings of the Links-Gould invariant}\label{SS:CablingsLG}

We recall some facts about cablings for details see \cite[\S 6.3.3]{BN}.  Let $L$ be a framed link.  For a non-zero integer $q$ let $L^{\otimes q}$ be the $q$th disconnected cabling of $L$.  If $f\in \V_{k}$ is a framed Vassiliev invariant of type $k$ then $f\circ(L\rightarrow L^{\otimes q})$ is a framed Vassiliev invariant of type $k$.  Let $V$ be a finite dimensional module of a Lie superalgebra $\g$ of type A-G.  Cablings of the invariant $Q_{\g,V}$ correspond to tensor powers of the module $V$, i.e.
\begin{equation}
\label{E:cablingsofgl}
W_{\g,V^{\otimes q},k}=[ Q_{\g,V,k}\circ(L\rightarrow L^{\otimes q}) ]
\end{equation}
where $[\cdot]$ is the map defined in subsection \ref{SS:KontT}.

Now we consider cablings of the Links-Gould invariant.  Let $\Valpha$ be the four dimensional simple module used to define $LG$.  Let $U$ be any finite dimensional $\gto$-module.  Let $Q^{\Valpha,U}$ be the invariant defined as follows.  If $T$ is a framed $(1,1)$-tangle then $T^{\otimes 2}$ is a (2,2)-tangle.  By closing the new parallel of $T^{\otimes 2}$ we have a framed $(1,1)$-tangle which we denote by $\cT$.  Label the original tangle of $\cT$ with $\Valphahat=\Valpha[[h]]$ and the new parallel with $U[[h]]$.  With this labeling on the $(1,1)$-tangle $\cT$ apply the R-T quantum group construction to obtain an element of $\End(\Valphahat)^{\gto}\cong \C[[h]] $.  We take this power series to be $Q^{\Valpha,U}(T)$.  We call $Q^{\Valpha,U}(T)$ the cabling of $LG$ and $Q_{\gto,U}$.

\begin{prop} Let $\stan$ is the standard $\gto$-module of column 
vectors and let $N_{\beta}$ be the 1-dimensional $\gto$ weight module with weight $(\beta,\beta |-\beta)$.
Let $U$ be any finite dimensional simple module of $\gto$.   
Then the invariant $Q_{\gto,V}$ is contained in the cablings of $LG$, $Q_{\gto,\stan}$ and $Q_{\gto,N_{\beta}}$.   In other words, any two knots distinguished by $Q_{\gto,V}$ are also distinguished by the cablings of $LG$, $Q_{\gto,\stan}$ and $Q_{N_{\beta}}$. 
\end{prop}

\begin{proof}
Throughout this proof we set $\g=\gto$. 
We will show that $W_{\g,U}$ can be constructed from cablings of the weight systems $\lgws$, $W_{\g,M}$ and $W_{\g,N_{\beta}}$.  Then Theorem \ref{T:gl}  and Theorem \ref{P:LGandWS} imply the desired result about invariants. 

First we define the weight system corresponding to $Q^{\Valpha, U}$.  
 Let $T$ is a framed $(1,1)$-tangle and let $D$ be a chord diagram on $T$.  By definition of the action of $\gto$ on the tensor product we have

\begin{equation}\label{D:ChordTensor}
\text{\includegraphics[width=2in, trim=0pt 0pt 0pt 0pt]{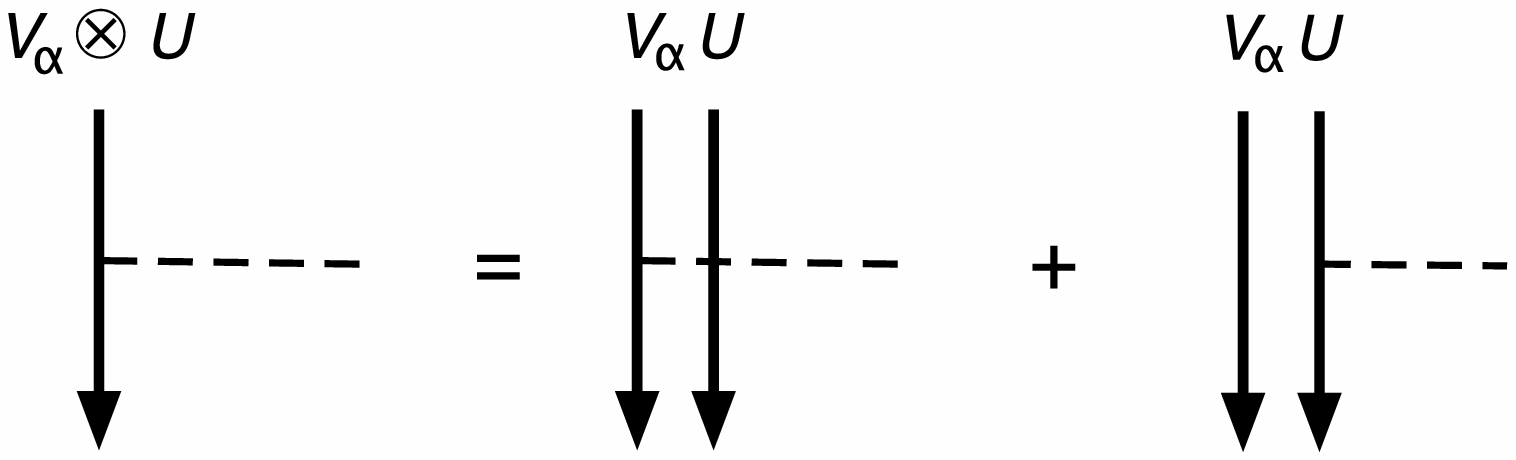}} 
\end{equation}

\noindent
It follows that the $G_{\g, \Valpha\otimes U}(D)$ can be constructed from a chord diagram on $T^{\otimes 2}$.  By closing the second parallel of $T^{\otimes 2}$ or equivalently contracting $U$ and $U^{*}$ we can associate $G_{\g, \Valpha\otimes U}(D)$ with an $\g$-invariant endomorphism of $\Valpha$ and thus a scalar.  Let $W^{\Valpha \otimes U}(D)$ be this scalar. 

We have the following extension of equality (\ref{E:cablingsofgl}) which holds for all $k\in \N$:
$$W^{\Valpha \otimes U}_{k}=[Q^{\Valpha,U}_{k}]$$
where $Q^{\Valpha,U}=\sum Q^{\Valpha,U}_{k}h^{k}$ and $ W^{\Valpha \otimes U}=(  W^{\Valpha \otimes U}_{k})_{k}$.  Therefore, $W^{\Valpha,U}$ is the weight system corresponding to the cabling of $LG$ and $Q_{\gto,U}$.

Next we will show that $W^{\Valpha \otimes U}$ is determined by weight systems coming from modules in a composition series of $\Valpha \otimes U$.  Denote the equivalence class of a module $L$ in the Grothendieck group by $[L]$.  Let $[U]=\sum_{\lambda\in I}[L_{\lambda}]$ where $L_{\lambda}$ is the simple module of weight $\lambda$ and $I$ is an index set.  
Recall that closing a colored $(1,1)$-tangle is the same as taking the supertrace of the corresponding endomorphism.  Also, taking the supertrace of an endomorphism only depends on the diagonal.  Therefore, we have 
\begin{equation}
\label{E:WisSupertrL}
W^{\Valpha\otimes U}(D)= \left( \sum_{\lambda\in I}supertrace(f_{L_{\lambda}})\right) c_{\Valpha}
\end{equation}
where $f_{M}\in \End(M)^{\g}$ and $c_{\Valpha}$ is a scalar corresponding to an endomorphism of $\Valpha$.  Equality (\ref{E:WisSupertrL}) says that $W^{\Valpha\otimes U}$ is determined by the cablings of weight systems arising from the simple modules $L_{\lambda}$ and the module $\Valpha$.  In other words, from (\ref{D:ChordTensor}) it follows that 
\begin{equation}
\label{E:WisTensorL}
W^{\Valpha\otimes U}=W^{\Valpha \bigotimes \underset{\lambda \in I}{\otimes}L_{\lambda}}.
\end{equation}  

Since $V_{0}$ is the trivial module we have $W^{V_{0} \otimes U}=W_{\gto,U}$.  This implies that the weight system $W_{\g, U}$ is contained in $W^{\Valpha\otimes U}$ and that it is determined by the set $\{L_{\lambda}: \lambda \in I\}$.  To finish the proof it suffices to show that the any simple module $L(\lambda)$ can be obtained in the Grothendieck group from a sum of the modules $\Valpha\otimes \stan^{k}\otimes N_{\beta}$.  The last statement follows from the following lemma:
%
%
\begin{lemma}\label{L:Groth}
    The Grothendieck group of finite dimensional $\gto$-modules is generated by the 
set of modules $A=\{  V_{\alpha} \otimes \stan^{\otimes k} 
\otimes N_{\beta} \: | \:  k\in \N , \alpha , \beta \in \C  \}.$
    \end{lemma}
\begin{proof}
We give a sketch of the proof of the lemma.  For a complete proof see Appendix of \cite{GThesis}. 

It suffices to show that every finite dimensional 
irreducible weight module $L(\gamma_{1},\gamma_{2}|\delta)$ of weight $(\gamma_{1},\gamma_{2}|\delta)$ is a linear combination 
of elements in $A$.  We will say a module 
is in the set $B$ if it is a linear combination of modules in A.

The proof follows from the observations:
\begin{enumerate}
  \item for any module $V$ in $B$ we have $V \otimes 
N_{\beta}$ is in $B$ for all $\beta \in \C$,  \label{I:Nbeta}
    \item for any $k \in \N$ and $\alpha \in \C$ the 
finite dimensional irreducible weight module of weight 
$(k,0|\alpha)$ is in $B$.\label{I:weightB}
\end{enumerate}
Observation (\ref{I:Nbeta}) allows us to shift the weight of a module.  In other words, if $V$ is a module of weight $(\gamma_{1},\gamma_{2}|\delta)$ then $V\otimes N_{\beta}$ is a module of weight $(\gamma_{1}+\beta,\gamma_{2}+\beta |\delta-\beta)$.  Therefore, if $B$ contains all the modules of weight  $(k,0|\alpha)$ the theorem follows.

Observation (\ref{I:weightB}) follows from strong induction on $k$.  Note that $\Valpha$ has highest weight $(0,0|\alpha)$ and so $\stan^{\otimes k}\otimes \Valpha$ has a highest weight $(k,0|\alpha)$.  Using the induction hypothesis one shows that $$[L(k,0|\alpha)]=[\stan^{\otimes k}\otimes V_{\alpha}]- \sum_{\lambda 
\in I}^{}  [L(\lambda)]  
$$
where $[V]$ is the equivalence class of a module $V$ in the Grothendieck group, $I$ is a particular index set and $[L(\lambda)]\in B$ for all $\lambda 
\in I$. 
\end{proof}    
In summary, Lemma \ref{L:Groth} and equality (\ref{E:WisTensorL}) imply that for every finite dimensional simple $\g$-module $L$ the weight system $W_{\g,L}$ is contained in the cablings of $W_{LG}$, $W_{\g,M}$ and $W_{\g,N_{\beta}}$.  Thus, Theorem \ref{T:gl}  and Theorem \ref{P:LGandWS} imply the desired result on the level of invariants.
\end{proof}

\subsection{The Links-Gould invariant and the colored HOMFLY}\label{SS:LGandHOMFLY}

It is known that the HOMFLY polynomial is built from $\mathfrak{gl}(n)$ and the standard module of column vectors (see \cite{RT}).  At the level of weight systems all knot invariants arising from finite dimensional $\mathfrak{gl}(n)$-modules are contained in the cablings of the HOMFLY polynomial.  Similarly, it is well known by experts \cite{Vain} that at the level of weight systems all the knot invariants arising from finite dimensional $\gl$-modules are contained in the cablings of the HOMFLY polynomial.  The following corollary is a direct consequence of Theorems \ref{T:gl}  and  \ref{P:LGandWS}.

\begin{corollary}\label{C:LGandHOMFLY}
Let $V$ be a finite dimensional $\gl$-module.  The invariant of knots $Q_{\gl,V}$ is contained in the colored HOMFLY.  In particular, every two knots distinguishable by $LG$ can be distinguished by the usual HOMFLY applied to some cabling of these knots.
\end{corollary}

\section{Two invariants arising from $\dtoa$}
Vogel \cite{Vogel} defined a Lie superalgebra $\dto$ over $R=\Q[a^{\pm},b^{\pm},c^{\pm}]/(a+b+c)$ which can be regarded as a generic version of the complex Lie superalgebra $\dtoa$ (by a change of of coefficients, one can recover $\dtoa$).  In \cite{Pat-Mir} Patureau-Mirand constructs to two link invariants $\Qbar$ and $\OZ$ arising from $\dto$.   These invariants are closely related to quantum invariants and weight systems composed with the Kontsevich integral, respectively.  
We will show that the results of \S \ref{SS:KontInv} imply that these invariants are the essentially the same.  This gives a positive answer to Conjecture 4.12 part 1 of \cite{Pat-Mir}.  Moreover, the invariant $\OZ$ is an invariant in three symmetric variables.  Our results imply that $\Qbar$ is an invariant which is symmetric in three variables.  This symmetry is not apparent at the level of the quantum group.  

We will now summarize Patureau-Mirand \cite{Pat-Mir} construction of the invariants $\Qbar$ and $\OZ$.  These construction are closely related to the constructions of $Q_{\dtoa, \adm}$ and $G_{\dtoa,\adm}\circ Z$, where $\adm$ is the adjoint module of $\dtoa$.  To define $\Qbar$ and $\OZ$ the normal constructions of the quantum group invariant and the weight systems arising from a module of a Lie superalgebra must be extended to the setting of 3-tangle (immersions of 1-3-valent graphs in $S^{3}$).  This involves accounting for trivalent vertices of the graph.  

Let $U_{h}(\dto)$ be the reparametrization of the quantization $U_{h}(\dtoa)$ given in Appendix of \cite{Pat-Mir}.  Let $\Qnet$ be the $\Q[e^{h}, e^{\alpha h}, e^{-h-\alpha h}]$-valued quantum group invariant of 3-tangles invariant arising from $U_{h}(\dto)$ and a irreducible $U_{h}(\dto)$-module $\dm$ whose classical limit is the adjoint module $\adm$ of $\dtoa$ (see \S 3.2 of \cite{Pat-Mir}).    In \cite{Vogel} Vogel define a functor $\Gnet$ from the category of chord diagrams on 3-tangles to the category of finite dimensional $\dto$-modules (also see \cite[\S 3.3]{Pat-Mir}).  Both $\Qnet$ and $\Gnet$ account for the trivalent vertices by mapping the diagram  \includegraphics[width=10pt]{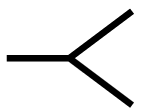} to the bilinear map $\dm\otimes \dm \rightarrow \dm$ (whose classical limit is the Lie bracket $\adm\otimes \adm\rightarrow \adm$).  
 
Let $Z_{ad}$ be the so-called Kontsevich-adjoint functor defined in \S 3.3.2 of \cite{Pat-Mir}.  The domain (range) of $Z_{ad}$ is the category of 3-tangle (resp. category of chord diagrams on 3-tangles).  Define $\ZGad=\Gnet \circ Z_{ad}$.  The restriction of $\ZGad$ to 3-nets (closed 3-tangles) is an invariant of 3-nets with values in the ring of symmetric series in three variables $\Q[[a_{1}, a_{2}, a_{3}]]^{\mathfrak{S}_{3}}$ modulo the ideal the ideal generated by $a_{1}+a_{2}+a_{3}$.  By specializing $a_{1}=h,  a_{2}=\alpha h$ and $a_{3}=-h -\alpha h$ 
we have $\ZGad$ is equal to $\Phi_{\dtoa} \circ Z_{ad}$ where $\Phi_{\dtoa}$ is a extension of $G_{\dtoa,\adm}$ to 3-tangles.  

We will now show that $\Qnet$ is equal to $\ZGad$.  After observing that trivalent vertices are preserved, this is essentially a restatement of Theorem \ref{T:KontInv} in the special case when the module $V$ is the adjoint module.  Setting $\g=\dtoa$ and $V=\adm$, let $E$ be the functor constructed in Theorem \ref{T:KontInv}.   By definition the functor $E$ preserves the morphism $\adm\otimes \adm\rightarrow \adm$ given by the bracket.  It follows that \begin{equation}
\label{E:ZadQtoa}
\Phi_{\dtoa} \circ Z_{ad}=Q^{3-tangles}_{\dtoa,\adm}.
\end{equation}  The invariant $\ZGad$ is determined by $\Phi_{\dtoa} \circ Z_{ad}$.  This is true because $\Phi_{\dtoa} \circ Z_{ad}$ is a specialization of $\ZGad$ such that each degree of the power series $\ZGad(T)$ is determined by infinitely many values of $\alpha$.  Thus, from (\ref{E:ZadQtoa}) we have that $\Qnet$ is equal to $\ZGad$. 

We will now give the definition of two 3-net invariants and use the last equality to show that the 3-net invariants are equal.  In \cite[\S4]{Pat-Mir} Patureau-Mirand showed that there exists two 3-net invariants $\OQ$ and $\OZ$ such that for any 3-net $L$ with $L={\displaystyle\ni}\circ T$ where $T$ is some 3-tangle we have
\begin{equation}
\label{E:OQOZ}
\scalebox{.8}{ \begin{picture}(0,0)%
\includegraphics{QZhalftTheta.pstex}%
\end{picture}%
\setlength{\unitlength}{3947sp}%
\begingroup\makeatletter\ifx\SetFigFont\undefined%
\gdef\SetFigFont#1#2#3#4#5{%
  \reset@font\fontsize{#1}{#2pt}%
  \fontfamily{#3}\fontseries{#4}\fontshape{#5}%
  \selectfont}%
\fi\endgroup%
\begin{picture}(3848,843)(826,-5320)
\put(4351,-4711){\makebox(0,0)[lb]{\smash{{\SetFigFont{14}{16.8}{\familydefault}{\mddefault}{\updefault}{\color[rgb]{0,0,0}$)$}%
}}}}
\put(3901,-5236){\makebox(0,0)[lb]{\smash{{\SetFigFont{14}{16.8}{\familydefault}{\mddefault}{\updefault}{\color[rgb]{0,0,0}$)$}%
}}}}
\put(826,-4711){\makebox(0,0)[lb]{\smash{{\SetFigFont{14}{16.8}{\familydefault}{\mddefault}{\updefault}{\color[rgb]{0,0,0}$\Qnet(T) =\OQ(L).\Qnet($}%
}}}}
\put(826,-5236){\makebox(0,0)[lb]{\smash{{\SetFigFont{14}{16.8}{\familydefault}{\mddefault}{\updefault}{\color[rgb]{0,0,0}$\ZGad(T) = \OZ(L).\Gnet($ }%
}}}}
\end{picture}%
}
\end{equation}
Since $\Qnet$ and  $\ZGad$ are equal it follows from  equations (\ref{E:OQOZ}) that $\OQ(L)=\OZ(L)$ for all 3-nets $L$.  This gives a positive answer to Part 1 of Conjecture 4.12 \cite{Pat-Mir}.  

We end this subsection by defining Patureau-Mirand's two link invariants and using the above results to show that they are related.  Patureau-Mirand \cite[\S4]{Pat-Mir} shows that $\OZ$ can be evaluated on links by the following equality:
\begin{equation*} \scalebox{.8}{ \begin{picture}(0,0)%
\includegraphics{ZIHX.pstex}%
\end{picture}%
\setlength{\unitlength}{3947sp}%
\begingroup\makeatletter\ifx\SetFigFont\undefined%
\gdef\SetFigFont#1#2#3#4#5{%
  \reset@font\fontsize{#1}{#2pt}%
  \fontfamily{#3}\fontseries{#4}\fontshape{#5}%
  \selectfont}%
\fi\endgroup%
\begin{picture}(4875,324)(526,-5323)
\put(526,-5236){\makebox(0,0)[lb]{\smash{{\SetFigFont{14}{16.8}{\familydefault}{\mddefault}{\updefault}{\color[rgb]{0,0,0}$\OZ\Big($}%
}}}}
\put(1876,-5236){\makebox(0,0)[lb]{\smash{{\SetFigFont{14}{16.8}{\familydefault}{\mddefault}{\updefault}{\color[rgb]{0,0,0}$\Big)=$}%
}}}}
\put(2401,-5236){\makebox(0,0)[lb]{\smash{{\SetFigFont{14}{16.8}{\familydefault}{\mddefault}{\updefault}{\color[rgb]{0,0,0}$\frac{1}{2}\OZ\Big($}%
}}}}
\put(5401,-5236){\makebox(0,0)[lb]{\smash{{\SetFigFont{14}{16.8}{\familydefault}{\mddefault}{\updefault}{\color[rgb]{0,0,0}$\big)\Big)$}%
}}}}
\put(4126,-5236){\makebox(0,0)[lb]{\smash{{\SetFigFont{14}{16.8}{\familydefault}{\mddefault}{\updefault}{\color[rgb]{0,0,0} $\frac{1}{2}\big($}%
}}}}
\end{picture}%
} \end{equation*}
Similarly for $\OQ$ he defines a link invariant $\Qbar$ by assigning a double point the following value:
\begin{equation*} \scalebox{.8}{ \begin{picture}(0,0)%
\includegraphics{QIHX.pstex}%
\end{picture}%
\setlength{\unitlength}{3947sp}%
\begingroup\makeatletter\ifx\SetFigFont\undefined%
\gdef\SetFigFont#1#2#3#4#5{%
  \reset@font\fontsize{#1}{#2pt}%
  \fontfamily{#3}\fontseries{#4}\fontshape{#5}%
  \selectfont}%
\fi\endgroup%
\begin{picture}(2850,324)(2551,-5323)
\put(5401,-5236){\makebox(0,0)[lb]{\smash{{\SetFigFont{14}{16.8}{\familydefault}{\mddefault}{\updefault}{\color[rgb]{0,0,0}$\big)\Big)$}%
}}}}
\put(4126,-5236){\makebox(0,0)[lb]{\smash{{\SetFigFont{14}{16.8}{\familydefault}{\mddefault}{\updefault}{\color[rgb]{0,0,0} $\frac{1}{2}\big($}%
}}}}
\put(2551,-5236){\makebox(0,0)[lb]{\smash{{\SetFigFont{14}{16.8}{\familydefault}{\mddefault}{\updefault}{\color[rgb]{0,0,0}$\OQ\Big($}%
}}}}
\end{picture}%
} \end{equation*}
Since $\OQ=\OZ$ it follows that for any link $L$ with $n$ components we have
$$\Qbar(L)=2\OZ(L) - n \OZ(O)$$
where $O$ is the unknot.  This gives a partial answer to Part 2 of Conjecture 4.12.

\Addresses\recd

\end{document}